\documentclass[a4paper]{amsart}
\usepackage[top=3.5cm, bottom=3.5cm, left=3cm, right=3cm]{geometry}
\usepackage{amsthm}
\usepackage{amsmath}
\usepackage[foot]{amsaddr}

\usepackage{amssymb}
\usepackage{mathtools}
\usepackage{graphicx} 
\usepackage[usenames]{color}
\usepackage{url,hyperref}
\usepackage{algorithm,algorithmic}
\usepackage{verbatim}
\usepackage{siunitx}
\usepackage{subfig}
\usepackage{flafter}
\usepackage{cite}
\usepackage{float}
\usepackage{multirow}

\usepackage{pgfplots}
\pgfplotsset{compat=1.16}
\usetikzlibrary{math}

\usepackage{academicons}
\newcommand{\orcid}[1]{\href{https://orcid.org/#1}{\textcolor[HTML]{A6CE39}{\aiOrcid}} #1}

\renewcommand{\orcid}[1]{\href{https://orcid.org/#1}{#1}}

\newenvironment{acknowledgements}
    {\section*{Acknowledgements}}
    {}
\newcommand{\norm}[1]{\|#1\|}
\usepackage{setspace}
\newtheorem{theorem}{Theorem}[section]
\newtheorem{lemma}[theorem]{Lemma}

\newtheorem{corollary}[theorem]{Corollary}

\theoremstyle{definition}
\newtheorem{definition}[theorem]{Definition}

\numberwithin{equation}{section}

\newcommand{\R}{\mathbb{R}}
\newcommand{\N}{\mathbb{N}}
\newcommand{\dx}{\,\mathrm{d}}

\newcommand{\python}{\textsc{Python}}
\newcommand{\fenics}{\textsc{FEniCS}}

\begin{document}

\title[Cond. recon. from power density data in lim. view]{Conductivity reconstruction from power density data in limited view}
\title[Conductivity from power density in limited view]{Conductivity reconstruction from power density data in limited view}

\author{Bjørn Jensen}
\address{
    Bjørn Jensen\\
    \orcid{0000-0002-4743-2631}\smallbreak
    Department of Mathematics and Statistics\\ 
    University of Helsinki\\ 
    00560 Helsinki, Finland
}
\email{bjorn.jensen@helsinki.fi}

\author{Kim Knudsen}
\address{
    Kim Knudsen\\
    \orcid{0000-0002-4875-3074}\smallbreak
    Department of Applied Mathematics and Computer Science\\
    Technical University of Denmark\\
    2800 Kgs. Lyngby, Denmark
}
\email{kiknu@dtu.dk}

\author{Hjørdis Schlüter}
\address{
    Hjørdis Schlüter\\ 
    \orcid{0000-0001-6659-3863}\smallbreak
    Department of Applied Mathematics and Computer Science\\
    Technical University of Denmark\\
    2800 Kgs. Lyngby, Denmark
}
\email{hjsc@dtu.dk}

\date{}

\begin{abstract}
In Acousto-Electric tomography, the objective is to extract information about the interior electrical conductivity in a physical body from knowledge of the interior power density data generated from prescribed boundary conditions for the governing elliptic partial differential equation.  In this note, we consider the problem when the controlled boundary conditions are applied only on a small subset of the full boundary. We demonstrate using the unique continuation principle that the Runge approximation property is valid also for this special case of limited view data. As a consequence, we guarantee the existence of finitely many boundary conditions such that the corresponding solutions locally satisfy a non-vanishing gradient condition. This condition is essential for conductivity reconstruction from power density data. In addition, we adapt  an existing reconstruction method intended for the full data situation to our setting. We implement the method numerically and investigate the opportunities and shortcomings when reconstructing from two fixed boundary conditions.\smallbreak

\textbf{Keywords:} {acousto-electric tomography, electrical impedance tomography, hybrid data tomography, coupled physics imaging, inverse problems, medical imaging}\smallbreak
\textbf{MSC2000:} {35R30; 65N21}
\end{abstract}

\maketitle

\section{Introduction}
Acousto-electric tomography (AET)~\cite{zhang2004a,Ammari_Bonnetier_Capdeboscq_Tanter_Fink_2008} is an imaging  modality for obtaining information about the interior electrical conductivity in a physical body. By combining electrostatic boundary measurements during an acoustic excitation of the body, we first recover the so-called interior power density data and from there reconstruct the conductivity. In principle, the modality yields high resolution and high contrast images, but the practical feasibility is yet to be demonstrated~\cite{jensen2019feasibility}.

The reconstruction problem in AET is fairly well understood in situations where the full boundary is available for the electrostatic measurements. In this manuscript we consider the more challenging problem in the limited view setting where we only have access to measurements on a possibly small part of the boundary. Our main result shows that reconstruction is feasible also in this setting, and in this way we generalize results pioneered in~\cite{bal2013inverse}.

We now define the involved quantities and the model linking them. The limited boundary control is modelled by a prescribed voltage on part of the domain boundary with a no-flux condition applied elsewhere. To this end we consider a connected open set $ \Omega \subset \R^d $, $ d = 2,3 $ with Lipschitz boundary $ \partial\Omega $ and let $ \Gamma \subset \partial\Omega $ be nonempty and open. We assume the conductivity denoted by $ \sigma $ satisfies $ \sigma \in L^\infty_+(\Omega) $, i.e. there exists $ \lambda \in (0,1) $ such that $ \lambda \leq \sigma \leq \lambda^{-1} $ a.e. in $ \Omega $. If $ d = 3 $ we assume further that $ \sigma $ is Lipschitz on $ \Omega $. The model is described by the partial differential equation
\begin{subequations} \label{eq:u}
    \begin{alignat}{2}
        Lu = -\nabla\cdot\sigma\nabla u &= 0 &\quad& \text{in $ \Omega $}, \label{eq:u-in-Omega} \\ 
        u &= f && \text{on $ \Gamma $}, \\
        L_\nu u = \sigma\partial_\nu u &= 0 && \text{on $ \partial\Omega\backslash\Gamma $,} \label{eq:u-neumann}
    \end{alignat}
\end{subequations}
where $ u $ is the electrical potential, $ \nabla u $ the electrical field and $ f \in H^\frac{1}{2}(\Gamma) := \{ v\vert_\Gamma : v \in H^1(\Omega) \} $ the prescribed voltage under our control. $ \nu $ denotes the outward unit normal to $ \partial\Omega $. From standard elliptic theory we know that~\eqref{eq:u} has a unique solution $ u \in H^1(\Omega) $ for each choice of boundary potential $ f $. 

In AET the data is given in terms of the mixed power density matrix $H \in \R^{d\times d}$ with elements
\begin{equation}
    H_{ij} = \sigma\nabla u_i\cdot\nabla u_j. \label{eq:pow-den}
\end{equation}
The function $ u_i $ denotes the solution to~\eqref{eq:u} for a boundary condition $ f_i $, $ i = 1,\dots, d. $
The inverse problem of AET is to uniquely and constructively recover the unknown conductivity $ \sigma $ from $H$. 

If we successfully pick $ f_i $ so that the electrical fields $ \nabla u_i $ are  linearly independent, then $ H $ satisfies the Jacobian constraint
\begin{equation} \label{eq:jac-constr}
    \det\Bigg[\begin{array}{ccc}\nabla u_1 & \cdots & \nabla u_d \end{array}\Bigg] \geq \delta > 0
\end{equation}
for some $ \delta > 0. $ This condition is crucial for the reconstruction problem in AET. 

There is a rich mathematical literature on AET for the full boundary case ($\Gamma = \partial\Omega$). For the two-dimensional problem, the Jacobian constraint is essentially guaranteed by the Rad\'o–Kneser–Choquet theorem using the two coordinate functions as boundary conditions. A reconstruction method is found in~\cite{bal2013inverse}. The anisotropic problem is considered in~\cite{BalMonard2012} with boundary conditions given in terms of complex geometrical optics (CGO) solutions generically depending on $\sigma$. In three dimensions, it is impossible to find three boundary conditions that can work for any conductivity~\cite{alberti2017critical}, but again a CGO construction can be used to provide four solutions~\cite{bal2013inverse}. A different approach using the Runge approximation property can be used for both the isotropic and anisotropic problem in any dimension~\cite{bal2014imaging} yielding a finite number of boundary conditions. Further work to establish smaller upper bounds on the number of boundary conditions can be found in~\cite{alberti2019combining}. We refer to~\cite{alberti2018lectures} for a complete overview. 

Our main result states that even in the limited view setting there are a finite number of boundary conditions such that locally the solutions satisfy ~\eqref{eq:jac-constr}. The method relies on the Runge approximation property that will be established for the mixed problem~\eqref{eq:u}. This approach is equivalent to the construction of so-called localized potentials~\cite{Gebauer2008}). We will in addition review the reconstruction method, emphasize the importance of~\eqref{eq:jac-constr}, and investigate the possibilities and shortcomings in the limited view setting.

In Section~\ref{sec:runge}, we recall the Runge approximation property and prove that the given PDE~\eqref{eq:u} satisfies the property. Section~\ref{sec:jac-constr} gives the solutions satisfying the Jacobian constraint. In Section~\ref{sec:recon} we review for completeness the reconstruction formula for $ \sigma $ in terms of  $H$ in the case $ d=2 $ following~\cite{bal2013inverse}.  We end the story in Section~\ref{sec:num} with a description of the numerical implementation and numerous computational experiments that shows the possibilities and limitations of limited view AET.

\section{The Runge approximation property} \label{sec:runge}

In this section, we recall the Runge approximation property and ascertain that the property holds for our PDE in question~\eqref{eq:u}. To prove that this property holds we make use of the unique continuation property, which is satisfied for the operators $ L $ and $ L_\nu $ as defined in~\eqref{eq:u}~\cite[Lemma 7.5]{alberti2018lectures}. Given an open non-empty subset $ \Sigma \subseteq \partial\Omega $ the unique continuation property states that if $ u \in H^1(\Omega) $ is a solution to $ Lu = 0 $ in $ \Omega $, which further satisfies $ u\vert_\Sigma = 0 $ and $ L_\nu u\vert_\Sigma = 0 $, then $ u $ vanishes in all of $ \Omega $; i.e. $ u = 0 $ everywhere in $ \Omega $.

We state the definition of the Runge approximation property.

\begin{definition}[Runge approximation property]\label{def:runge-apprx}
    We say that $ L $ satisfies the \emph{Runge approximation property} if for any simply connected Lipschitz domain $ \Omega' \Subset \Omega $
    and any $ u \in H^1(\Omega') $ such that $ Lu = 0 $ in $ \Omega' $ there exists a sequence $ u_n \in H^1(\Omega) $ such that (a) $ Lu_n = 0 $ in $ \Omega $, and (b) $ u_n\vert_{\Omega'} \to u $ in $ L^2(\Omega') $.
\end{definition}

It follows from the satisfaction of the unique continuation property that our PDE satisfies the Runge approximation property as well. The following lemma is a modified version of a similar result in~\cite[Thm. 7.7]{alberti2018lectures} and the proof follows the same general lines as well:

\begin{theorem} \label{thm:L-runge-approx-weak}
Let $ \Omega \subset \R^d $ be a Lipschitz bounded domain and $ L $ as defined in~\eqref{eq:u}. Then $ L $ satisfies the Runge approximation property with the sequences $ (u_n) $ satisfying the boundary restriction  $ L_\nu u_n = 0 $ on $ \partial\Omega\backslash\Gamma $ for all $ n \in \N $.
\end{theorem}

\begin{proof}
Assume without loss of generality that $ \Omega $ is connected. Take $ \Omega' \Subset \Omega $ as in Definition~\ref{def:runge-apprx} and $ u \in H^1(\Omega') $ such that
$$ Lu = 0 \quad \text{in $ \Omega' $}. $$
Set $ F := \{ v\vert_{\Omega'} : v \in H^1(\Omega), Lv = 0\ \text{in $ \Omega $}, L_\nu v = 0\ \text{on $ \partial\Omega\backslash\Gamma $}\} $. Suppose by contradiction that the Runge approximation property does not hold that we managed to pick a $ u $ for which there is no sequence in $ F $ converging to $ u $ in $ L^2 $-norm. By the Hahn-Banach theorem there exists a functional $ g \in L^2(\Omega')^\ast $ such that $ g(u) \neq 0 $ and $ g(v) = 0 $ for $ v \in F $. In other words, there exists $ g \in L^2(\Omega') $ such that $ (g,u)_{L^2(\Omega')} \neq 0 $ and $ (g,v)_{L^2(\Omega')} = 0 $ for all $ v \in F $.

Consider now the extension by zero of $ g $ to $ \Omega $, which by an abuse of notation is still denoted by $ g $. Let $ V := \{ \phi\vert_\Gamma : \phi \in H^\frac12(\partial\Omega)\} $ and fix $ \phi \in V $. Let $ w,v \in H^1(\Omega) $ be the unique solutions to
\begin{equation*}
    \begin{aligned}
        Lw &= g && \text{in $ \Omega $}, \\
        w &= 0 && \text{on $ \Gamma $}, \\
        L_\nu w &= 0 && \text{on $ \partial\Omega\backslash\Gamma $},
    \end{aligned} \quad\quad \text{and} \quad\quad 
    \begin{aligned}
        Lv &= 0 && \text{in $ \Omega $}, \\
        v &= \phi && \text{on $ \Gamma $}, \\
        L_\nu v &= 0 && \text{on $ \partial\Omega\backslash\Gamma $}.
    \end{aligned}
\end{equation*}
By definition of $ g $ there holds $ (g,v)_{L^2(\Omega)} = 0 $. Thus, integration by parts shows
\begin{equation*}
    0 = -(v,g)_{L^2(\Omega)} = (Lv,w)_{L^2(\Omega)} - (v,Lw)_{L^2(\Omega)} = \int_\Gamma (L_\nu w)\phi\dx\sigma.
\end{equation*}
Since the above holds for all $ \phi \in H^\frac{1}{2}(\Gamma) $ we obtain $ L_\nu w = 0 $ on $ \Gamma $. Observe now that $ w $ is solution to $ Lw = 0 $ in $ \Omega\backslash\Omega' $ such that $ w = 0 $ and $ L_\nu w = 0 $ on $ \Gamma $. In view of the unique continuation property 
we have $ w = 0 $ in $ \Omega\backslash\Omega' $, therefore $ w = 0 $ and $ L_\nu w = 0 $ on $ \partial\Omega' $. As a result, by integrating by parts we obtain
\begin{align*}
    \int_{\Omega'} gu\dx{x} &= \int_{\Omega'} (Lw)u\dx{x}\\ &= \int_{\Omega'} \underbrace{(Lu)}_{=0}w\dx{x} + \int_{\partial\Omega'} \underbrace{(L_\nu w)}_{=0}u\dx{s} + \int_{\partial\Omega'} (L_\nu u)\underbrace{w}_{=0}\dx{s}\\ &= 0,
\end{align*}
where the last identity follows from the definition of $ u $. This contradicts the assumptions on $ g $, since $ (g,u)_{L^2(\Omega')} \neq 0 $.
\end{proof}

By the application of Schauder estimates, Sobolev embedding and standard elliptic theory, assuming Lipschitz regularity for the coefficient yields the following result.
\begin{lemma} \label{lem:appx-lem}
    Assume $ \sigma $ is Lipschitz and let $ \Omega' \Subset \Omega $, $ x_0 \in \overline{\Omega'} $ and $ s \in (0,\operatorname{dist}(\Omega',\partial\Omega)) $. Let $ u_0 \in C^{1,\alpha}(\Omega) $, $ \alpha \in (0,1)$ satisfy
    \[
        -\Delta u_0 = 0, \quad \text{in $ B(x_0,s) $}.
    \]
    
    Then for any $ \delta > 0 $ there is an $ r \in (0,s) $ and $ u_{x_0,\delta} \in H^1(\Omega) $ satisfying $ Lu_{x_0,\delta} = 0 $ in $ \Omega $ and $ L_\nu u_{x_0,\delta} = 0 $ on $ \partial\Omega\backslash\Gamma $ such that
    \[
        \|u_{x_0,\delta} - u_0\|_{C^1(\overline{B(x_0,r)})} \leq \delta.
    \]
\end{lemma}
The proof can be found in Appendix A. As any linear map $ u_0 $ satisfies the conditions of the lemma, the corollary is immediate.
\begin{corollary} \label{cor:vec-approx}
    Let $ \sigma $, $ \Omega' $ and $ x_0 $ be as in Lemma~\ref{lem:appx-lem}. For any $ \delta > 0 $ and $ \mathbf v \in \R^d $ there is $ r \in (0,1) $ with $ B(x_0,r) \subset \Omega' $ and $ u_{x_0,\delta} \in H^1(\Omega) $ satisfying $ Lu_{x_0,\delta} = 0 $ in $ \Omega $ and $ L_\nu u_{x_0,\delta} = 0 $ on $ \partial\Omega\backslash\Gamma $ such that
    \[
        \|\nabla u_{x_0,\delta} - \mathbf v\|_{C^0(B(x_0,r))} \leq \delta.
    \]
\end{corollary}

\section{Satisfying the Jacobian constraint} \label{sec:jac-constr}
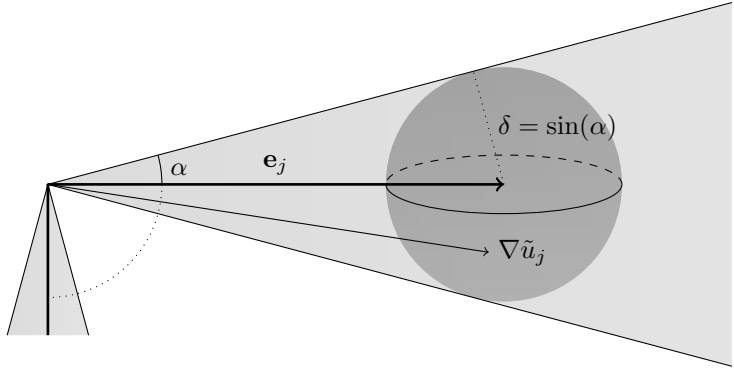
\begin{figure}[!ht]
    \centering
    \begin{tikzpicture}
\fill[
  top color=gray!50,
  bottom color=gray!10,
  shading=axis,
  opacity=0.35
  ] 
  (0,0) circle (1.55291427cm);
\fill[
  left color=gray!50,
  right color=white,
  shading=axis,
  opacity=0.125
  ] (3, -2.41154273188) -- (-6,0) -- (3, 2.41154273188);
\draw
  (-1.55291427,0) arc (180:360:1.55291427cm and 0.38822856765cm);
\draw[dashed]
  (1.55291427,0) arc (0:180:1.55291427cm and 0.38822856765cm);
\draw[->,line width=1] (-6,0) -- node[above] {$\mathbf{e}_j $} (0,0);
\draw[->] (-6,0) -- (-.2,-.9) node[right] {$\nabla\tilde u_j $};

\draw (3,-2.41154273188) -- (-6,0) -- (3, 2.41154273188);
\draw[] (-4.5,0) arc (0:15:1.5cm) node[midway, anchor=west] {$ \alpha $};
\draw[dotted] (0,0) -- ++(-0.40192378864668, 1.5) node[midway, anchor=west] {$ \delta = \sin(\alpha)$};

\fill[
  left color=gray!50,
  right color=white,
  shading=axis,
  opacity=0.125
  ] (-6-0.53589838486,-2) -- (-6,0) -- (-6+0.53589838486, -2);
\draw (-6-0.53589838486,-2) -- (-6,0) -- (-6+0.53589838486, -2);
\draw[line width=1] (-6,0) -- (-6,-2);
\draw[dotted] (-4.5,0) arc (-0:-90:1.5cm);
\end{tikzpicture}
    \caption{Illustration of the cone $ C_j $ and the ball of approximating vectors at most $ \delta = \sin(\alpha) $ away from $ \mathbf e_j $
    .}
    \label{fig:cone}
\end{figure}

We summarize in this section how the above guarantees solutions that satisfies the Jacobian constraint. The treatment follows along the lines of~\cite[Sec. 7.3]{alberti2018lectures}:

First note that the canonical basis vectors $ \mathbf e_{j} $ satisfies  $\det\begin{bmatrix} \mathbf e_1 \dots \mathbf e_d \end{bmatrix} = 1, $ i.e.\ the Jacobian constraint. Next, fix a sufficiently small $\delta>0$. Corollary~\ref{cor:vec-approx} provides the existence of solutions to~\eqref{eq:u} (including suitable boundary conditions) that can approximate  $v=\mathbf e_{j}$ locally in small balls. This is illustrated in Figure~\ref{fig:cone} if considered from the same relative origin then $ \nabla \tilde u_j $ is contained in a cone $ C_j. $ Consequently, the solutions satisfy locally the Jacobian constraint. Finally, the union  of such small balls is an open cover of $ \Omega' $ and can, due to compactness of $\overline{\Omega'}$, be exhausted to a finite open cover with $M\in \mathbb{N} $ balls. As each ball requires a set of $ d $ boundary conditions to satisfy,  we ultimately end up with $ Md $ boundary conditions.


The theorem is as follows:
\begin{theorem}\label{thm:jac-constr-satisf}
Assume $ \sigma $ is Lipschitz continuous in $ \overline{\Omega}$, and consider any compactly embedded domain $ \Omega' \Subset \Omega. $ Then there exists a finite set of boundary conditions $ \{\phi_j\}_{j=1}^{Md} $ with corresponding solutions $ \{u_j\}_{j=1}^{Md} $ to~\eqref{eq:u} with $ f = \phi_{j} $,  such that at any point $ x_0 \in {\Omega'} $ there is an open neighborhood $ V_{x_0} $ of $ x_0 $ and a subset $ \{u_{n_1},\dots, u_{n_d}\} $ satisfying
\[ \det\Bigg[\begin{array}{ccc}\nabla u_{n_1}(x) & \cdots & \nabla u_{n_d}(x) \end{array}\Bigg] \geq \frac{1}{2}, \quad \text{for all $ x \in V_{x_0}$} \]
\end{theorem}
\begin{proof}
Let $ \delta > 0 $ be arbitrary. For each $ x_0 \in \Omega', $ and $v =  \mathbf e_j, $  Corollary~\ref{cor:vec-approx} gives an $ r_j > 0 $ and an approximating solution $ \tilde u_j \in H^1(\Omega) $ satisfying $ L\tilde u_j = 0 $ in $ \Omega $ and $ L_\nu\tilde u_j = 0 $ on $ \partial\Omega\backslash\Gamma $ such that $ \|\nabla \tilde u_j - \mathbf e_j\|_{C^0(B(x_0,r_j))} < \delta $ in $B(x_0,r_0)$ with $r_0 = \min(r_1,r_2,\ldots,r_d).$  By continuity of the determinant, by fixing $ \delta $ sufficiently small we can guarantee
$$ \left|\det\begin{bmatrix} \nabla\tilde u_1 \dots \nabla\tilde u_d \end{bmatrix} - \det\begin{bmatrix} \mathbf e_1 \dots \mathbf e_d \end{bmatrix}\right| < \frac{1}{2}, $$ and by the reverse triangle inequality
$$ \left|\det\begin{bmatrix} \nabla\tilde u_1 \dots \nabla\tilde u_d \end{bmatrix} \right| > \frac{1}{2}. $$ 
Since $\Omega'\Subset \Omega,$ there is a finite subcover $\{B(x_m,r_m)\}_{m=1}^M$ of $\Omega'$ that yields the $Md$ boundary conditions.
\end{proof}

\section{Reconstruction formula} \label{sec:recon}
In this section we review without proofs the main formulas used for the reconstruction procedure based on~\cite{monard2012a} in order to prepare for the numerical experiments in section \ref{sec:num}. We consider the case $ d = 2 $; for $ d=3 $ we refer to~\cite{bal2013inverse}. Throughout this section we consider a power density matrix $ \mathbf H = (H_{ij})_{1\leq i,j\leq 2} $ for which the corresponding solutions $u_1$ and $u_2$ are assumed to satisfy the Jacobian constraint~\eqref{eq:jac-constr}.

The derivation goes in 2 steps. The aim in the first step is to obtain $\mathbf J_i = \sqrt{\sigma}\nabla u_i $, $ i = 1,2 $, from $ H_{ij} $. $\mathbf J_i $ is the current density up to a factor of $ \sqrt\sigma $; i.e. the current density is precisely $ \sqrt\sigma \mathbf J_i $. In a slight abuse of the name, we call the first step the ``current density step''. The second step is the recovery of an equation for $ \sigma $ from $\mathbf J_i. $ We call this the ``conductivity step''.

\subsection{Current density step} 
To this end consider the matrices $ \mathbf J = \begin{bmatrix} \mathbf J_1 & \mathbf J_2 \end{bmatrix} $ and $ \mathbf H = (H_{ij})_{1\leq i,j\leq 2} $; the latter obviously symmetric. Under the assumed positive lower bound on $ \sigma $ and the Jacobian constraint~\eqref{eq:jac-constr} it is easy to see that for any $ \mathbf x = (x_1,x_2) \neq 0 $, $ \mathbf x^T\mathbf H\mathbf x = \sigma|\mathbf q_1 + \mathbf q_2|^2 > 0 $, where $ \mathbf q_i = x_i\nabla u_i \neq 0 $ and $ \mathbf q_1\neq -\mathbf q_2 $. Hence $ \mathbf H $ is positive definite.

The matrix $\mathbf J$ is then orthonormalized into a $SO(2)$-valued matrix $\mathbf R $ via a transformation of the form $\mathbf R = \mathbf J \mathbf T^T$. By definition $\mathbf R$ is orthogonal and $\det \mathbf R=1$.
As $ \mathbf R $ is a rotation matrix, it is parameterized by
\begin{equation*}
    \mathbf R(\theta) = \begin{bmatrix} \cos\theta & -\sin\theta \\ \sin\theta & \cos\theta \end{bmatrix},
\end{equation*}
dependent only on a single parameter $ \theta $. This definition implies that when the exact matrices $\mathbf T$ and $\mathbf J$ are known the function $\theta$ can be computed by
\begin{equation}\label{eq:theta}
    \theta=\text{argument}(R_{11}+i R_{21}).
\end{equation}
Defining $ \mathbf T = (T_{ij})_{1\leq i,j\leq 2} $ and $ \mathbf T^{-1} = (T^{ij})_{1\leq i,j\leq 2} $, and letting 
\[ 
    \mathbf{V_{ij}} = \nabla(T_{i1})T^{1j}+\nabla(T_{i2})T^{2j},
\]
then $ \theta $ is determined by the following equation~\cite[Eq. (65)]{monard2012a}:
\begin{equation}
    \nabla\theta = \mathbf F, \label{eq:gradTheta}
\end{equation} 
with 
\[
    \mathbf F=\frac{1}{2}\left(\mathbf{V_{12}} -  \mathbf{V_{21}} - \mathcal{J}\nabla\log D\right),
\]
$\mathcal J = \begin{bmatrix} 0 & -1 \\ 1 & 0 \end{bmatrix},$ and  $D= (H_{11}H_{22}-H_{12}^2)^{\frac{1}{2}}.$  Applying the divergence operator to \eqref{eq:gradTheta} yields the Poisson equation
\begin{equation}\label{thetaPois}
    \left\{\begin{aligned}
            &\Delta \theta = \nabla \cdot \mathbf F \,\, \text{in }\Omega,\\
            &\theta \text{ given on }\partial \Omega. 
    \end{aligned}\right.
\end{equation}
We discuss in section~\ref{sec:Tandtheta} below that we can assume $\theta$ to be known on the whole boundary $ \partial\Omega. $

\subsection{Conductivity step} 
Reconstruction of $\sigma$ is based on~\cite[eq. (68)]{monard2012a}
\begin{align} 
    &\hspace{6mm}\nabla \log \sigma = \mathbf{G} \label{eq:grad-sigma}
\end{align}
with 
\begin{align*} 
\mathbf G &=\cos(2\theta) \mathbf K + \sin(2 \theta)\mathcal{J}\mathbf K\\
\mathbf K &=\mathcal{U}(\mathbf{V_{11}}-\mathbf{V_{22}})+ \mathcal{J}\mathcal{U}(\mathbf{V_{12}}+\mathbf{V_{21}}) \quad \text{and} \quad \mathcal{U}=\begin{bmatrix} 1 & 0\\ 0 & -1 \end{bmatrix}.\nonumber
\end{align*}
 We assume that $\sigma$ is known on the whole boundary $\partial \Omega$ and  then reconstruct $\sigma$ by solving the Poisson equation
\begin{equation}\label{sigPois}
    \left\{\begin{aligned}
        &\Delta (\log(\sigma)) = \nabla \cdot \mathbf G \,\, \text{in }\Omega,\\
        &\log(\sigma) \text{ given on }\partial \Omega. 
    \end{aligned}\right.
\end{equation}

\subsection{Choice of \texorpdfstring{$\mathbf{T}$}{T} and knowledge of \texorpdfstring{$\theta$}{theta}}\label{sec:Tandtheta}
The decomposition $\mathbf R= \mathbf J \mathbf T^T$ is not unique as any transfer matrix $\mathbf T_1$ that satisfies this equation for some rotation matrix $\mathbf R_1$ gives rise to another rotated transfer matrix $\mathbf T_2 = \mathbf{\overline{R}}^T \mathbf T_1$ that satisfies the same equation for the corresponding rotation matrix $\mathbf R_2=\mathbf R_1 \mathbf{\overline{R}}$. This holds true for any rotation matrix $\mathbf{\overline{R}}$. In theory, the choice of $\mathbf T$ will not influence the reconstruction procedure, as every choice of $\mathbf{T}$ with corresponding $\theta$ will work to extract the functionals $\mathbf{J}_i=\sqrt{\sigma} \nabla u_i$ from the entries of $H_{ij}=\sigma \nabla u_i \cdot \nabla u_j$. However, numerically a simple choice of $\mathbf{T}$ can be an advantage. For this reason, we choose Gram-Schmidt orthonormalization to obtain the following $\mathbf T$, as in this case the functionals $\mathbf{V_{ij}}$ have the simplified form as in \eqref{eq:Vexp}:
\begin{equation}\label{eq:Tdef}    
    \mathbf T = \begin{bmatrix} H_{11}^{-\frac{1}{2}} & 0 \\ -H_{12} H_{11}^{-\frac{1}{2}} D^{-1} &H_{11}^{\frac{1}{2}} D^{-1} \end{bmatrix}.
\end{equation} 
By the Jacobian constraint~\eqref{eq:jac-constr}, $H_{11}>0$ and thus $\mathbf T$ is well-defined. For this choice of $\mathbf{T}$ the function $\theta$ is given by the angle between $\nabla u_1$ and the $x_1$-axis, as in this case the first column of $\mathbf R$ simplifies to
\begin{equation*}
    \mathbf{R_1}=T_{11}\mathbf{J_1} + T_{12}\mathbf{J_2} = \dfrac{\nabla u_1}{\vert\nabla u_1\vert}.
\end{equation*}
Using \eqref{eq:theta} we then calculate
\begin{align*}
    \theta=\text{argument}(\partial_1 u_1 + i \partial_2 u_1).
\end{align*}
In addition, the vector fields $\mathbf{V_{ij}}$ can be written explicitly in terms of $\mathbf{H}$:
\begin{align}\label{eq:Vexp}
\begin{aligned}
    \mathbf{V_{11}}&= \nabla \log H_{11}^{-\frac{1}{2}}, & \mathbf{V_{12}}&=0,\\ \mathbf{V_{21}}&= -\frac{H_{11}}{D} \nabla \left( \frac{H_{12}}{H_{11}}\right), & \mathbf{V_{22}}&=\nabla \log \left( \frac{H_{11}^{\frac{1}{2}}}{D}\right).
    \end{aligned}
\end{align}
An outline of the reconstruction procedure using this choice for $\mathbf{T}$ is shown in Algorithm~\ref{algo}.

Knowledge of $\theta$ at the boundary is  essential  in the reconstruction procedure. The gradient equation~\eqref{eq:gradTheta} can be solved in two different ways for $\theta$: Either by knowledge of $\theta$ at a single point and integrating along curves originating from that point, or by solving the Poisson problem~\eqref{thetaPois}, where knowledge of $\theta$ is required on the whole boundary. As we suggest the second option, the question is whether knowledge of $\theta$ at the boundary is a valid assumption. Information about $\theta$ is related to the gradient $\nabla u_1$ and hence the current $\sigma \nabla u_1$, as both vector fields have the same direction. The functional $\sigma \nabla u_1$ can be decomposed into two parts with contribution from the unit outward normal $\nu$ and its direct orthogonal $\eta=\mathcal{J}\nu$:
\begin{equation*}
        \sigma\nabla u_1= (\sigma\nabla u_1 \cdot \nu)\nu + (\sigma\nabla u_1 \cdot \eta)\eta.
\end{equation*}
Along the part of the boundary $\partial \Omega \backslash \Gamma$ the Neumann data is known, so that knowledge of $(\sigma\nabla u_1 \cdot \eta)\eta$ is required to determine $\sigma\nabla u_1$ along this part of the boundary. Along $\Gamma$ the Dirichlet data and thus the directional derivative $(\nabla u_1 \cdot \eta)$ is known, so that knowledge of $\sigma$ and $\sigma\nabla u_1 \cdot \nu$ is required to determine $\sigma\nabla u_1$ along $\Gamma$. This requires the following information on the different parts of the boundary:
\begin{itemize}
    \item Knowledge of $\sigma$ and the potential $u_1$ yield information about $\theta$ on $\partial \Omega\backslash \Gamma$
    \item Knowledge of $\sigma$ and $\sigma\nabla u_1 \cdot \nu$ yield information about $\theta$ on $\Gamma$
\end{itemize}
However, along $\partial \Omega \backslash \Gamma$ the vanishing Neumann data gives us the additional information that $\sigma\nabla u_1$ is solely determined by the direct orthogonal $\eta$. 
Hence, on $\partial \Omega \backslash \Gamma$ the direction of $\sigma\nabla u_1$ is either $\eta$ or $-\eta$. The sign of $\eta$ is influenced by the choice of the boundary function $f_1$. By theorem \ref{thm:Max} (the weak maximum principle) the maximum and minimum of $u_1$ is attained at points $x_{\text{max}}$ and $x_{\text{min}}$ on the Dirichlet boundary $\Gamma$. At these points $\nabla u_1$ will then point outwards and inwards of the domain respectively, as the gradient points in the direction of highest increase. The maxima and minima of $f_1$ will therefore influence the direction of the current and this information can be used to determine sign of $\eta$ and thus $\theta$ on $\partial \Omega \backslash \Gamma$.

\begin{algorithm}[H]
\begin{flushleft} 
Generating data: Choose a set of boundary conditions $(f_1,f_2)$ so that the corresponding solutions satisfy the Jacobian constraint \eqref{eq:jac-constr}. Use the solutions to generate the power density matrix $\mathbf H$.
\end{flushleft}
\begin{enumerate}
        \item Use the measurement matrix $\mathbf H$ to define $\mathbf T$ and the vector fields $\mathbf{V_{ij}}$ as in \eqref{eq:Tdef} and \eqref{eq:Vexp} respectively
        \item Reconstruct $\theta$ by solving the boundary value problem \eqref{thetaPois}
        \item Reconstruct $\sigma$ by solving the boundary value problem \eqref{sigPois}
\end{enumerate}
 \caption{Reconstruction procedure}\label{algo}
\end{algorithm}

\section{Numerical Examples} \label{sec:num}
The \textsc{Matlab} and \textsc{Python} code to generate the numerical examples can be found on \textsc{GitLab}: \href{https://lab.compute.dtu.dk/hjsc/conductivity-reconstruction-from-power-density-data-in-limited-view.git}{GitLab code}.\par
\begin{figure}[!ht]
    \centering
    \begin{minipage}[t]{0.5\textwidth}
        \centering
        \includegraphics[width=\linewidth]{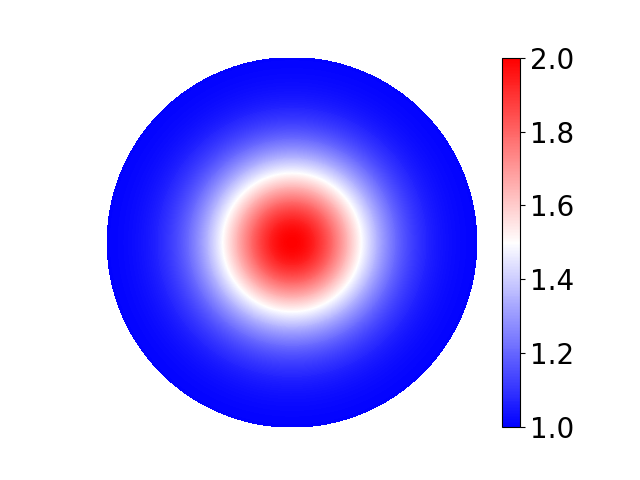}
        \caption*{Case 1}
    \end{minipage}%
    \begin{minipage}[t]{0.5\textwidth}
        \centering
        \includegraphics[width=\linewidth]{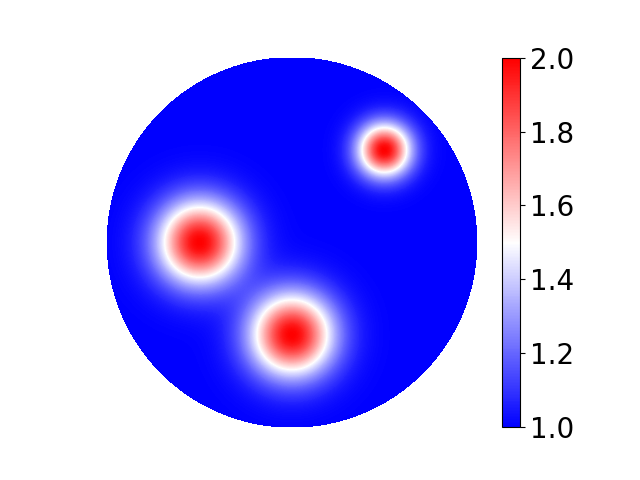}
        \caption*{Case 2}
    \end{minipage}
    \caption{The conductivities $\sigma$ used for the reconstruction procedure.}
    \label{fig:TrueSig}
\end{figure}

To investigate the performance of the reconstruction algorithm in the limited view setting, we implement the algorithm in \python{} and used \fenics{}~\cite{fenics} to solve the PDEs. Here we mainly focus on performance for two boundary conditions. We use a fine mesh to generate our power density data and a coarser mesh to address the reconstruction problem. Unless otherwise stated we used $N_{\text{data}}=44880$ nodes in the high-resolution case, while for the smaller mesh we considered a resolution of $N_{\text{recon}}=20100$ nodes. For both meshes, we use $\mathbb P_1 $ elements. We consider the domain $\Omega$ to be the unit disk: $\Omega=B(\mathbf 0, 1)$. Furthermore, we consider two test cases for the conductivity $\sigma$ defined by:
\begin{align*}
    \sigma_{\text{case 1}}(x^1,x^2) &= 1 + e^{-5 \left( (x^1)^2+(x^2)^2\right)},\\
    \sigma_{\text{case 2}}(x^1,x^2) &= 1 + e^{-20 \left( \left(x^1+\frac{1}{2}\right)^2+(x^2)^2\right)} + e^{-20 \left( \left(x^1\right)^2+\left(x^2+\frac{1}{2}\right)^2\right)}\\
    &\hspace{6mm}+ e^{-50 \left( \left(x^1- \frac{1}{2}\right)^2+\left(x^2-\frac{1}{2}\right)^2\right)},
\end{align*}
for $(x^1,x^2)\in \Omega$. Figure~\ref{fig:TrueSig} illustrates the conductivities.\par
To test the reconstruction procedure we follow the procedure as outlined in Section~\ref{sec:recon}. To generate the power density data we simply consider the coordinate functions as boundary conditions:
\begin{equation*}
    (f_1,f_2)=(x^1,x^2), \quad (x^1,x^2) \in \Gamma.
\end{equation*}
We make this choice, as there is no theory developed on how to constructively choose the boundary functions in the case of mixed Dirichlet and Neumann conditions such that \eqref{eq:jac-constr} is satisfied. However, for Dirichlet boundary conditions and when the domain is convex and 2-dimensional it is known that the coordinate functions $(f_1,f_2)=(x_1,x_2)$ yield solutions $(u_1,u_2)$ so that these satisfy the Jacobian condition (see e.g. \cite{BaumanMariniNesi01}). Motivated by this theory we limit ourselves to these boundary conditions.\par
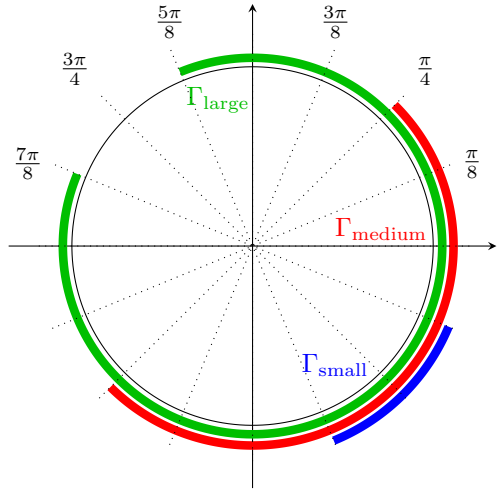
\begin{figure}[!ht]
    \centering
    \begin{minipage}[t]{0.5\textwidth}
        \centering
        \begin{tikzpicture}
\begin{axis}[
	trig format plots=rad,
	axis x line=center,
    axis y line=center,
    xlabel={},
    ylabel={},
    ticks=none,
    unit vector ratio*=1 1 1,
    xmin=-1.35, xmax=1.35,
    ymin=-1.35, ymax=1.35,
    height=8cm,
]
\tikzmath{\v=pi/8;} 

\addplot[domain=0:2*pi, samples=200] ({cos(x)},{sin(x)});

\addplot[domain=7*pi/8:2*pi+5*pi/8, samples=200, green!75!black, line width=.75ex] ({1.05*cos(x)},{1.05*sin(x)}) node [pos=.175,xshift=7ex,yshift=24ex,anchor=north west] {$ \Gamma_{\textup{large}}$};
\addplot[domain=2*pi-3*pi/4:2*pi+pi/4, samples=200, red, line width=.75ex] ({1.1125*cos(x)},{1.1125*sin(x)}) node [pos=.5,xshift=4ex, yshift=16ex,anchor=north east] {$ \Gamma_{\textup{medium}}$};
\addplot[domain=2*pi-3*pi/8:2*pi-1*pi/8, samples=200, blue, line width=.75ex] ({1.175*cos(x)},{1.175*sin(x)}) node [pos=.5,xshift=-1.8ex,yshift=.5ex,anchor=south east] {$ \Gamma_{\textup{small}}$};

\addplot[domain=-1.2:1.2, dotted] ({cos(\v)*x},{sin(\v)*x}) node [anchor=west] {$ \frac{\pi}{8}$};
\addplot[domain=-1.2:1.2, dotted] ({cos(2*\v)*x},{sin(2*\v)*x}) node [anchor=south west] {$ \frac{\pi}{4}$};
\addplot[domain=-1.2:1.2, dotted] ({cos(3*\v)*x},{sin(3*\v)*x}) node [anchor=south] {$ \frac{3\pi}{8}$};
\addplot[domain=-1.2:1.2, dotted] ({cos(4*\v)*x},{sin(4*\v)*x});
\addplot[domain=-1.2:1.2, dotted] ({cos(5*\v)*x},{sin(5*\v)*x}) node [anchor=south] {$ \frac{5\pi}{8}$};
\addplot[domain=-1.2:1.2, dotted] ({cos(6*\v)*x},{sin(6*\v)*x}) node [anchor=south east] {$ \frac{3\pi}{4}$};
\addplot[domain=-1.2:1.2, dotted] ({cos(7*\v)*x},{sin(7*\v)*x}) node [anchor=east] {$ \frac{7\pi}{8}$};
\addplot[domain=-1.2:1.2, dotted] ({cos(8*\v)*x},{sin(8*\v)*x});

\end{axis}
\end{tikzpicture}
    \end{minipage}
    \caption{Different sizes of $\Gamma$ used for the reconstruction procedure.}
    \label{fig:GammaSizeLoc}
\end{figure}

  
We consider different sizes for the subset of the boundary, $\Gamma$, that we can control, and show these in Figure~\ref{fig:GammaSizeLoc}. For each choice of $\Gamma$ we then solve the two PDEs using \fenics{} and use the computed solutions $u_1$ and $u_2$ to generate the 3 power densities $H_{11}, H_{12}$ and $H_{22}$. Furthermore, for comparison we compute the true angle $\theta$ from knowledge of the true gradient $\nabla u_1$. \par

\subsection{Reconstruction of \texorpdfstring{$\theta$}{theta}} 

\begin{figure}[!ht]
    \centering
    \begin{minipage}[t]{0.5\textwidth}
        \centering
        \includegraphics[width=\linewidth]{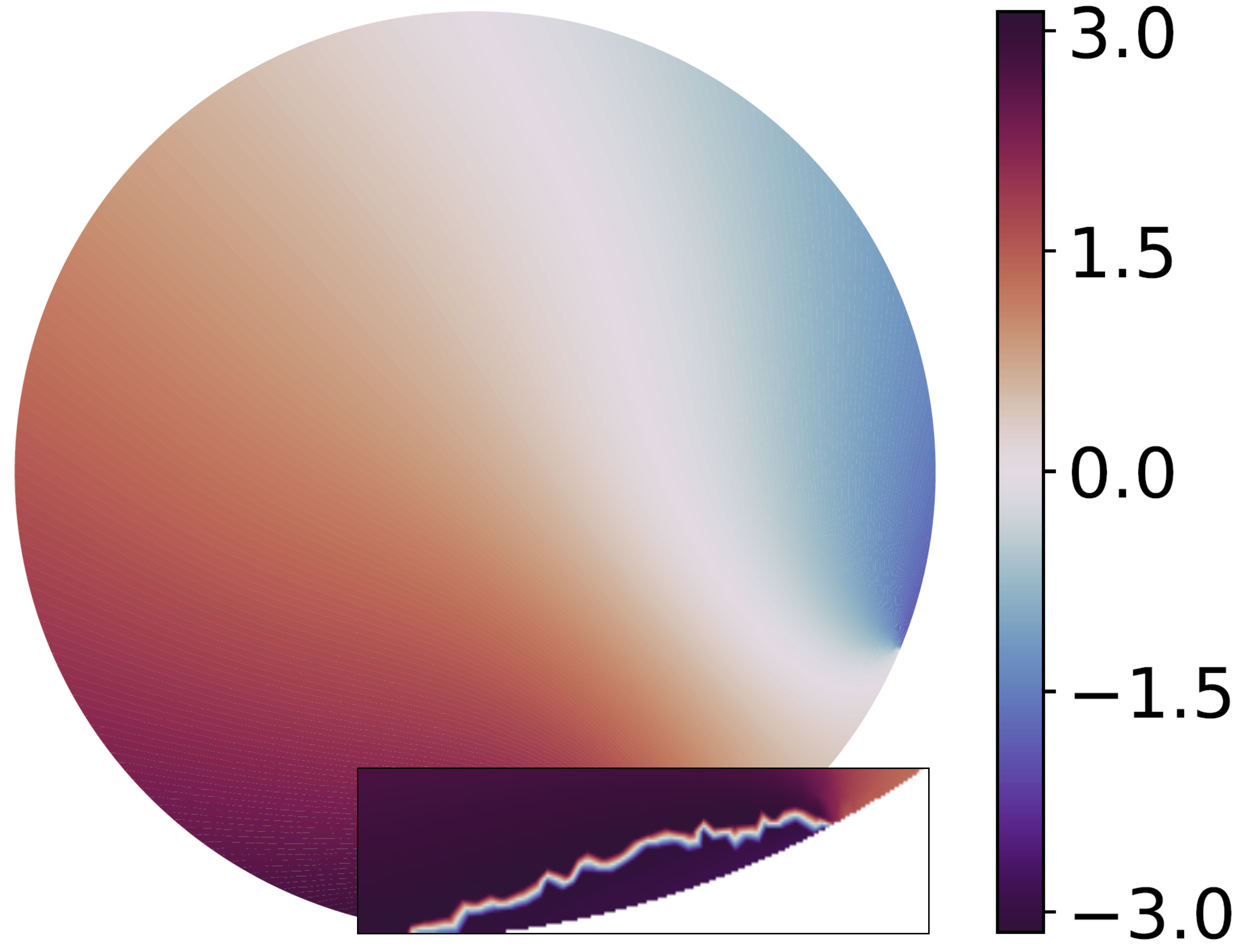}
        \caption*{(a) True $\theta$ \vspace{4mm}}
    \end{minipage}%
    \begin{minipage}[t]{0.5\textwidth}
        \centering
        \includegraphics[width=\linewidth,trim={0 3mm 0 0},clip]{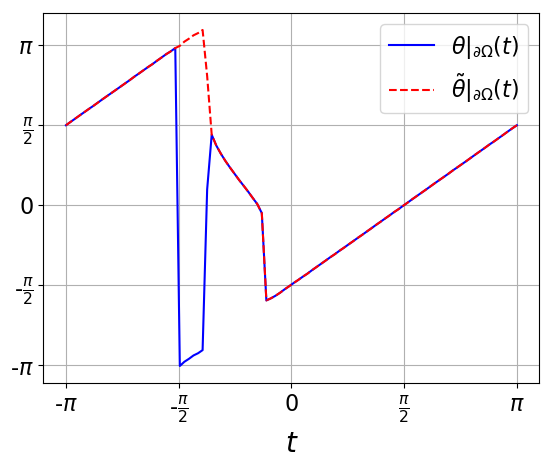}
        \caption*{(b) True $\theta$ and a modified version $\tilde{\theta}$ along the boundary $\partial\Omega$.}
    \end{minipage}
    \caption{True expression for $\theta$ assigned to a smooth function space for $\sigma$ as in test case 1 and having control over $\Gamma_{\text{small}}$. The right figure compares the same function to a modified version $\tilde{\theta}$ defined in \eqref{eq:modThet} along the boundary.} 
    \label{fig:Theta}
\end{figure}

For the three different sizes of $\Gamma$, $\Gamma_{\text{large}}, \Gamma_{\text{medium}}$ and $\Gamma_{\text{small}}$, as seen in Figure~\ref{fig:GammaSizeLoc}, we reconstruct $\theta$ for $\sigma$ defined as in test case 1. As $\theta$ is a function mapping from $ \Omega$ to the unit circle $ S^1 $, which is a periodic space, we choose on $ \partial\Omega $ a representation of $\theta$ which forms a continuous curve when unwrapping $S^1$ to an interval of $\R$. We denote this representation $\tilde{\theta}$ and it is constructed as follows.
\begin{equation}\label{eq:modThet}
    \tilde{\theta}=\begin{cases}\theta + 2\pi & t\in \left[-\frac{\pi}{2},-\frac{3\pi}{8} \right]\\
    \theta & \text{otherwise}.
    \end{cases}
\end{equation}
The representation is illustrated in Figure \ref{fig:Theta}(a). This is done to avoid any discontinuous transitions, which could otherwise occur from $\theta$ taking values on alternating sides of the end points of the interval; this problem is illustrated in Figure \ref{fig:Theta}(b). The special region in which we make the modification to obtain continuity is identified by considering the boundary conditions. As $ \theta $ represents the angle between $-\eta$ and the $x_1$-axis a discontinuity, when viewing $S^1$ in $\R$ can be identified to occur at exactly the points $ t = -\frac\pi2 $ and $ t - \frac{3\pi}8 $. 

\begin{figure}[!htb]
    \centering
    \begin{minipage}[t]{0.33\textwidth}
        \centering
        \includegraphics[width=\linewidth]{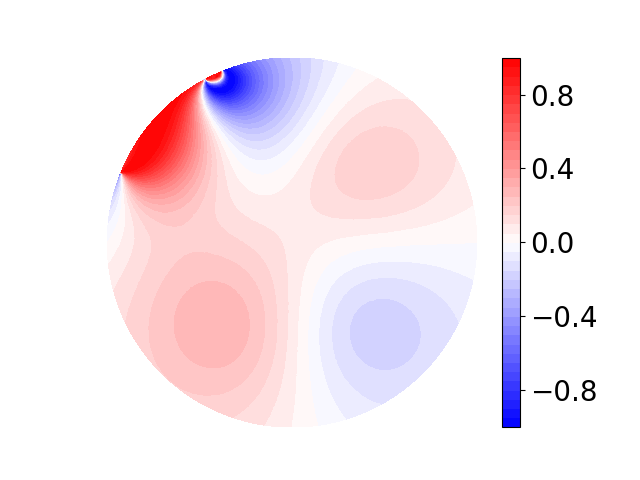}
    \end{minipage}%
    \begin{minipage}[t]{0.33\textwidth}
        \centering
        \includegraphics[width=\linewidth]{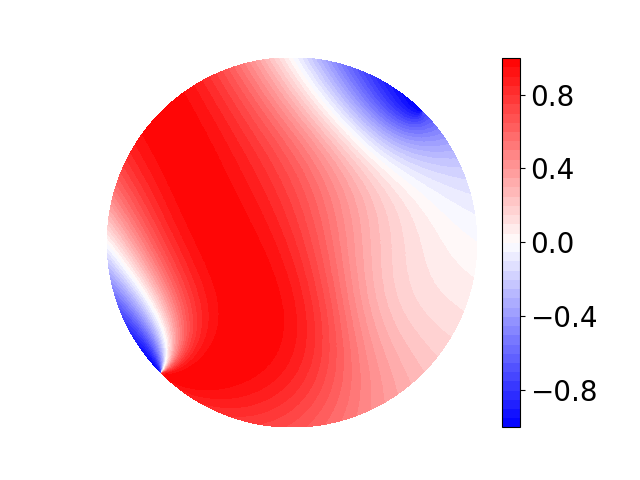}
    \end{minipage}%
    \begin{minipage}[t]{0.33\textwidth}
        \centering
        \includegraphics[width=\linewidth]{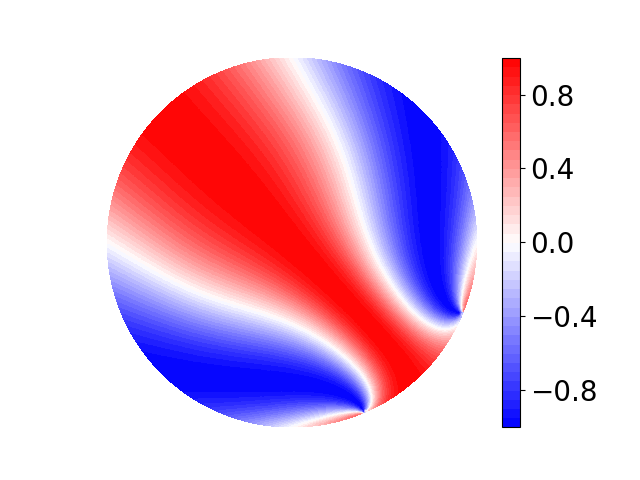}
    \end{minipage}
    \begin{minipage}[t]{0.33\textwidth}
        \centering
        \includegraphics[width=\linewidth]{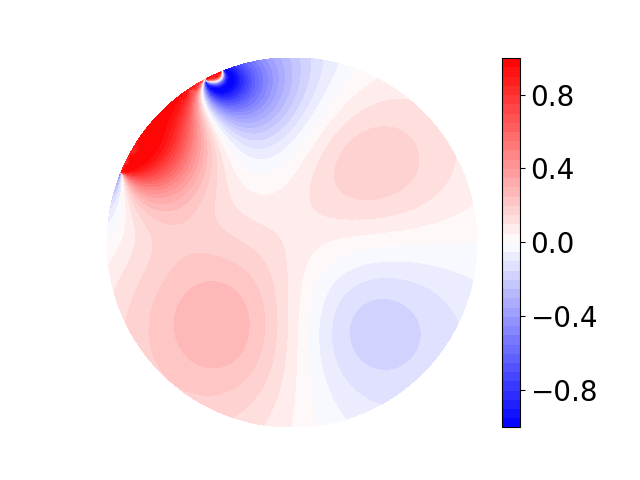}
        \caption*{$\Gamma_{\text{large}}$}
    \end{minipage}%
    \begin{minipage}[t]{0.33\textwidth}
        \centering
        \includegraphics[width=\linewidth]{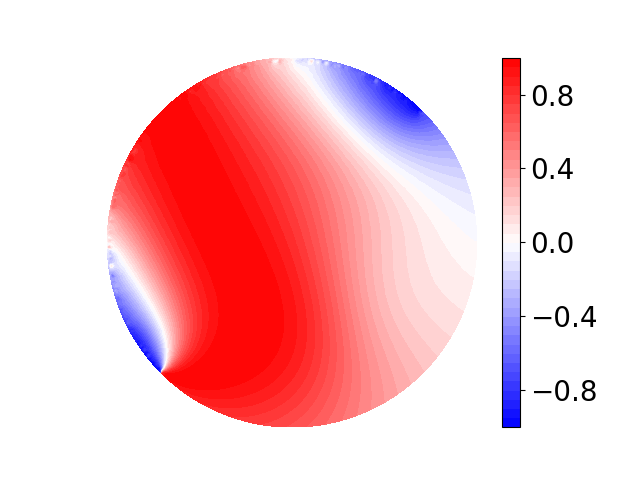}
        \caption*{$\Gamma_{\text{medium}}$}
    \end{minipage}%
    \begin{minipage}[t]{0.33\textwidth}
        \centering
        \includegraphics[width=\linewidth]{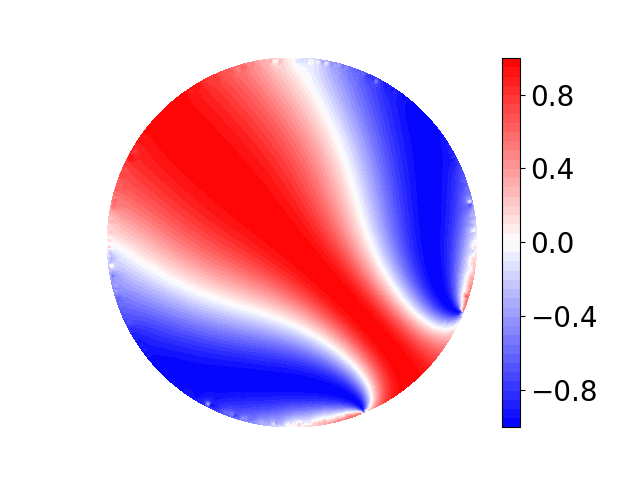}
        \caption*{$\Gamma_{\text{small}}$}
    \end{minipage}
    \caption{True expression for $\sin(2\theta)$ in the upper row and reconstructions in the lower row for varying sizes of $\Gamma$ and $\sigma$ being as in test case 1.}
    \label{fig:thetaRecs}
\end{figure}

Using the representation $\tilde{\theta}$ as a boundary condition, we reconstruct $\theta$ for the three different sizes of $\Gamma$. 
The reconstructions of $ \theta $ are compared to the true expressions in Figure~\ref{fig:thetaRecs}, while table~\ref{tab:RelErrSize} shows the relative $L^2$-error. For the smaller sizes of $\Gamma$ there appear artifacts close to the boundary $\partial \Omega \backslash \Gamma$ otherwise there is no visual distinction between true expressions and reconstructions, which is reflected in the relative errors as well. Note that due to the above modifications we use $\sin(2\theta)$ instead of $\theta$ for comparison.

\begin{table}[!ht]
    \centering
    \caption{Relative $L^2$ errors for varying sizes of $\Gamma$.}
    \begin{tabular}{p{23mm} p{9mm}||p{22mm}|p{22mm}|p{22mm}}
        && $\Gamma_{\text{large}}$ & $\Gamma_{\text{medium}}$ & $\Gamma_{\text{small}}$\\ \hline \hline
        \multirow{2}{*}{Min det$(\mathbf{H})$} & case 1 & $3.94 \cdot 10^{-6}$ & $3.87 \cdot 10^{-10}$ & $9.94 \cdot 10^{-18}$\\
        & case 2 & $2.94 \cdot 10^{-6}$ & $3.57 \cdot 10^{-10}$ & $1.07 \cdot 10^{-17}$\\\hline
        Rel. $L^2$ error & case 1 & (0.79\%,2.04\%) & (1.40\%,2.01\%) & (2.24\%,2.37\%)\\
        $(\cos(2\theta),\sin(2\theta))$ & case 2 & (0.77\%,1.86\%) & (1.41\%,1.97\%) & (2.25\%,2.33\%)\\ \hline
        \multirow{2}{*}{Rel. $L^2$ error $\sigma$} & case 1 & 32.02\% & 104\% & 177\%\\
        &case 2 & 33.62\% & 108\% & 180\%\\\hline
    \end{tabular}
    \label{tab:RelErrSize}
\end{table}

\subsection{Reconstruction of \texorpdfstring{$\sigma$}{sigma}}
\begin{figure}[!ht]
    \centering
    \begin{minipage}[t]{0.33\textwidth}
        \centering
        \includegraphics[width=\linewidth]{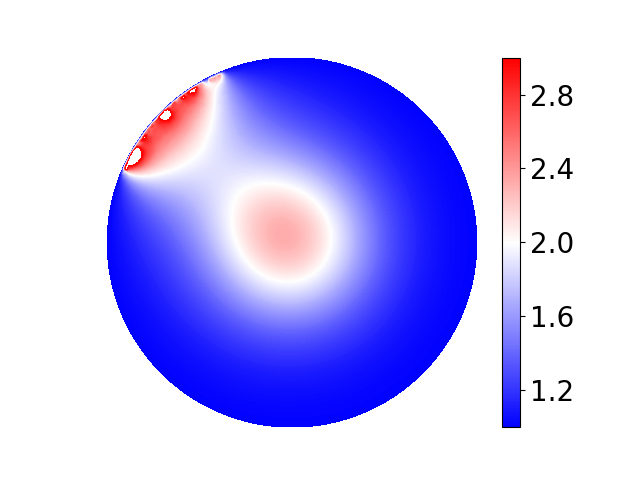}
    \end{minipage}%
    \begin{minipage}[t]{0.33\textwidth}
        \centering
        \includegraphics[width=\linewidth]{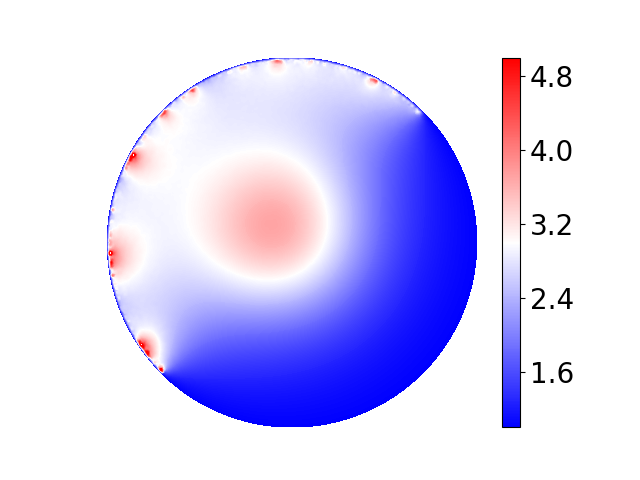}
    \end{minipage}%
    \begin{minipage}[t]{0.33\textwidth}
        \centering
        \includegraphics[width=\linewidth]{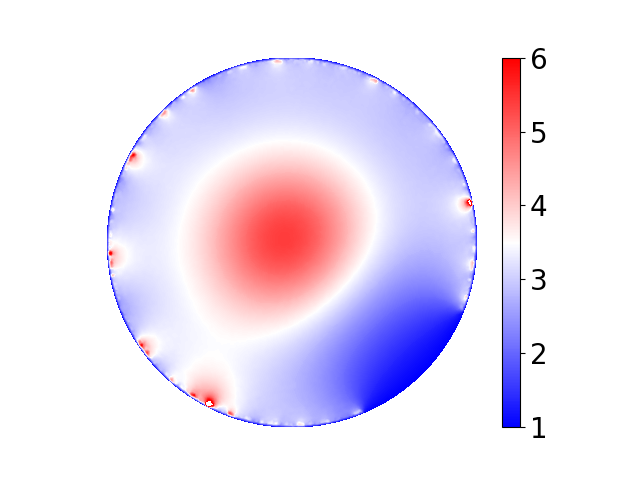}
    \end{minipage}
    \begin{minipage}[t]{0.33\textwidth}
        \centering
        \includegraphics[width=\linewidth]{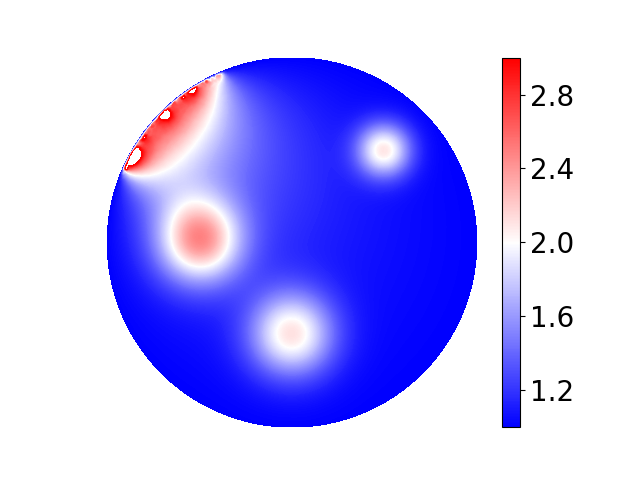}
        \caption*{$\Gamma_{\text{large}}$}
    \end{minipage}%
    \begin{minipage}[t]{0.33\textwidth}
        \centering
        \includegraphics[width=\linewidth]{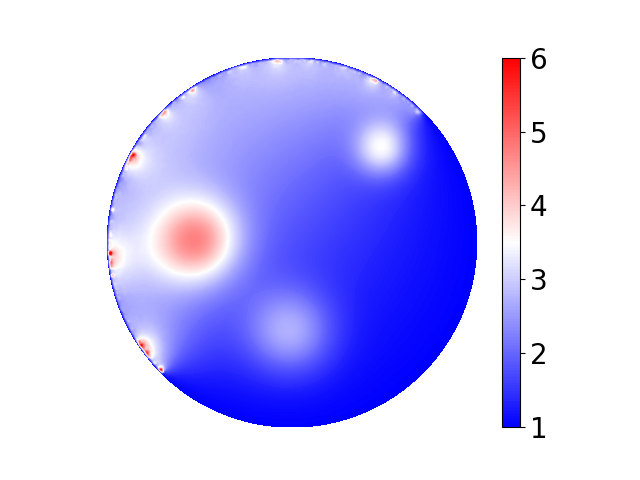}
        \caption*{$\Gamma_{\text{medium}}$}
    \end{minipage}%
    \begin{minipage}[t]{0.33\textwidth}
        \centering
        \includegraphics[width=\linewidth]{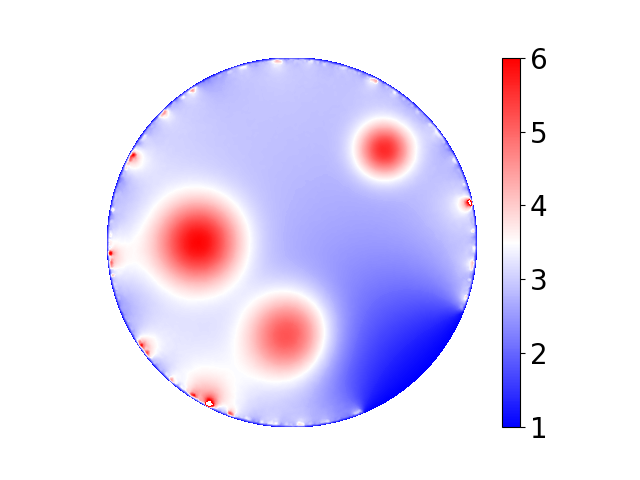}
        \caption*{$\Gamma_{\text{small}}$}
    \end{minipage}
    \caption{Reconstructions of $\sigma$ as in test case 1 for varying sizes of $\Gamma$ in the upper row and for test case 2 in the lower row. Note that different color scales are used.}
    \label{fig:sigmaRecs}
\end{figure}

\begin{figure}[!ht]
    \centering
    \begin{minipage}[t]{0.33\textwidth}
        \centering
        \includegraphics[width=\linewidth]{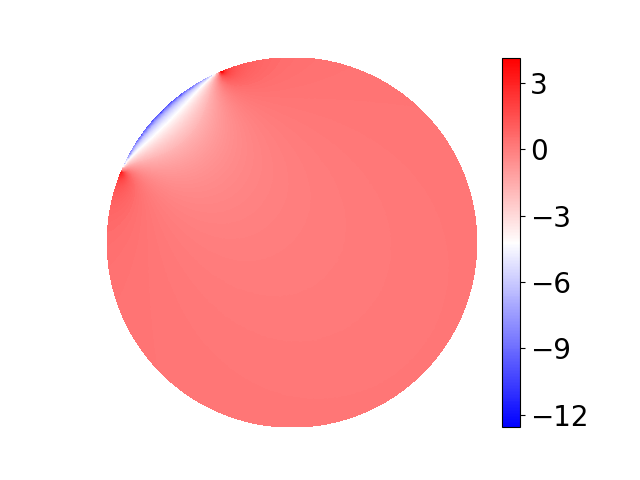}
        \caption*{$\Gamma_{\text{large}}$}
    \end{minipage}%
    \begin{minipage}[t]{0.33\textwidth}
        \centering
        \includegraphics[width=\linewidth]{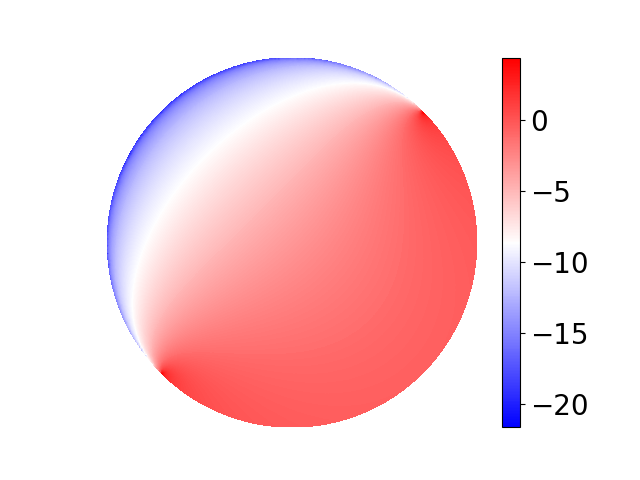}
        \caption*{$\Gamma_{\text{medium}}$}
    \end{minipage}%
    \begin{minipage}[t]{0.33\textwidth}
        \centering
        \includegraphics[width=\linewidth]{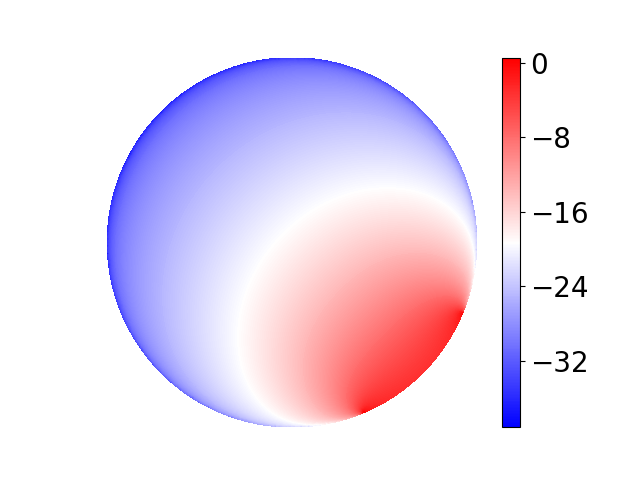}
        \caption*{$\Gamma_{\text{small}}$}
    \end{minipage}
    \caption{The expression $\log(\det \mathbf{H})$ for varying sizes of $\Gamma$. Large negative values (blue regions) correspond to values of $\det \mathbf{H}$ close to zero. }
    \label{fig:Dets}
\end{figure}

We reconstruct $\sigma$ for both test cases in Figure~\ref{fig:TrueSig} and for the three different sizes of $\Gamma$. Here we use the reconstructed $\theta$ from the previous step. The results are shown in Figure~\ref{fig:sigmaRecs}; note that the intervals for the colorbars are different, as for smaller sizes of $\Gamma$ there is more variation in the reconstructed $\sigma$.  Furthermore, Table~\ref{tab:RelErrSize} shows the relative $L^2$-error in the reconstructions. Both visually and quantitatively, we see that the quality of the reconstructed $\sigma$ decrease for decreasing size of $\Gamma$. Comparing the different test cases of $\sigma$, from the relative error we see that the quality of the reconstruction is slightly better for test case 1. The slight difference in quality may be explained by the fact that in test case 2 there are smaller features closer to the boundary. Overall, the decreasing quality for decreasing size of $\Gamma$ may follow from the fact that the assumption of a non-vanishing determinant of $H$ is almost violated for $\Gamma_{\text{small}}$. Figure~\ref{fig:Dets} illustrates the expression $\log(\det(\mathbf{H}))$ for all sizes of $\Gamma$ and their minimal values can be found in Table~\ref{tab:RelErrSize}. We see from the Figure that with decreasing size of the boundary of control, the values of $\det(\mathbf{H})$ decrease close to $\partial \Omega \backslash \Gamma$. By assumption, the determinant should be non-vanishing everywhere in the domain, but for $\Gamma_{\text{small}}$ the minimal value is $9.94 \cdot 10^{-18},$ and in comparison to values of the determinant elsewhere in the domain, this is effectively zero. As a consequence, the reconstruction in such regions of low determinants are not reliable, even though the algorithm does produce a result (see reconstructions in Figure \ref{fig:sigmaRecs} below.) 
\par

\begin{figure}
    \centering
\includegraphics[width=0.8\textwidth]{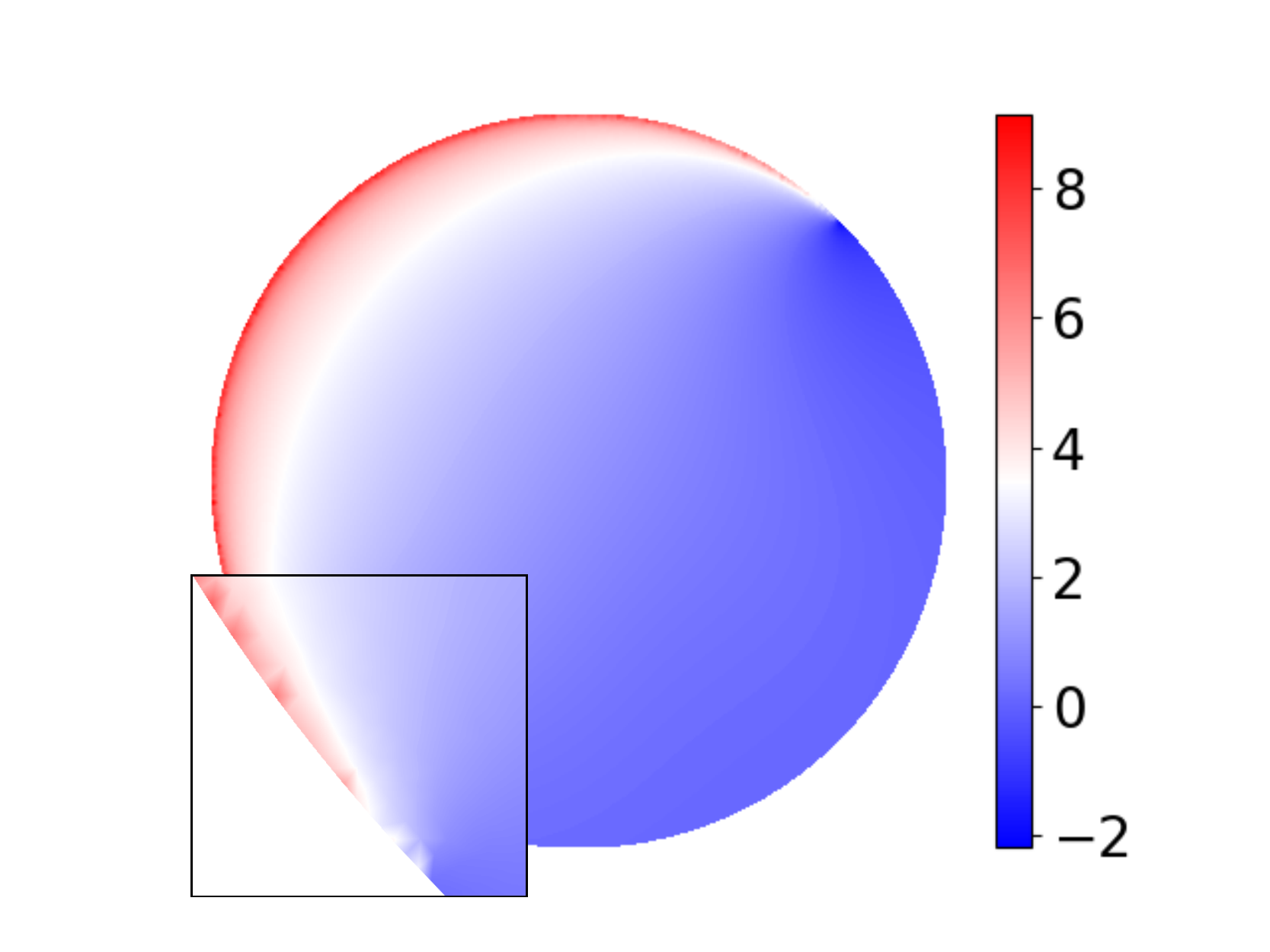}
    \caption{The fraction $\log \left( \frac{\sqrt{H_{11}}}{D} \right)$ in test case 1 when having control over $\Gamma_{\text{medium}}$. We zoom in on some of the discontinuities towards the boundary $\partial \Omega \backslash \Gamma$.}
    \label{fig:Frac}
\end{figure}

\begin{figure}[!ht]
    \centering
        \begin{minipage}[t]{0.33\textwidth}
        \centering
        \includegraphics[width=\linewidth]{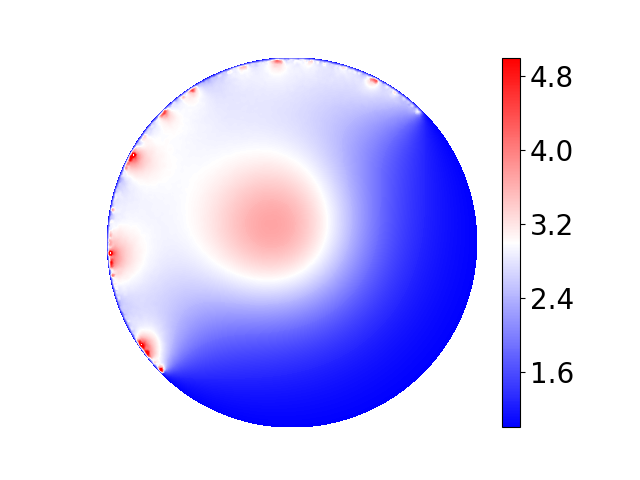}
        \caption*{$N_{\text{small}}$}
    \end{minipage}%
    \begin{minipage}[t]{0.33\textwidth}
        \centering
        \includegraphics[width=\linewidth]{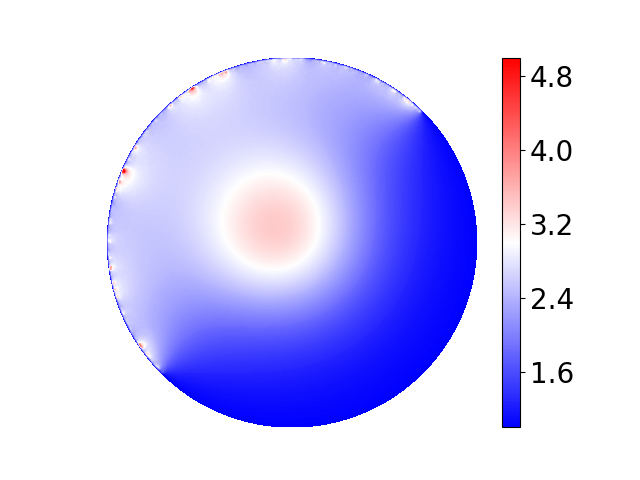}
        \caption*{$N_{\text{medium}}$}
    \end{minipage}%
    \begin{minipage}[t]{0.33\textwidth}
        \centering
        \includegraphics[width=\linewidth]{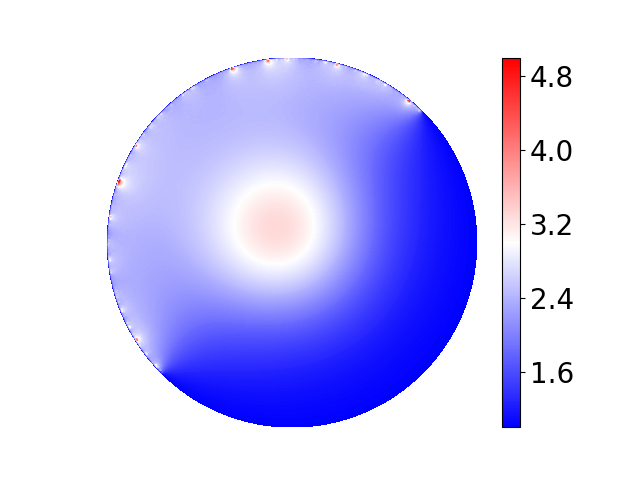}
        \caption*{$N_{\text{large}}$}
    \end{minipage}
    \caption{Reconstructions of $\sigma$ as in test case 1 for varying mesh sizes when having control over $\Gamma_{\text{medium}}$. }
    \label{fig:mesh}
\end{figure}

The size of $\Gamma$ roughly reflects how large an amount of the domain can be controlled: While for $\Gamma_{\text{large}}$ the reconstruction looks reasonable for most of the domain, for $\Gamma_{\text{medium}}$ this is the case for about half of the domain and for $\Gamma_{\text{small}}$ this does only apply for a small part of the domain. For the last two cases the part of the domain that is difficult to control is dominated by values of the reconstructed conductivity that are way higher than the true values and this is reflected in the high relative errors above $100\%$. What is common for all sizes of $\Gamma$ is that there appear artifacts close to the part of the boundary that cannot be controlled, $\partial \Omega \backslash \Gamma$. For $\Gamma_{\text{medium}}$ and $\Gamma_{\text{small}}$ the artifacts occur at distinct points along the boundary. Further experiments showed that the locations of the artifacts seem not influenced by the choice of the boundary conditions. This indicates that the appearance of the artifacts is not caused by the PDE. However, we think that they are a mesh-dependent phenomenon. To understand what part in the reconstruction procedure induces the artifacts, we illustrate the fraction $\log \left( \frac{\sqrt{H_{11}}}{D} \right)$ in Figure~\ref{fig:Frac}. This fraction appears in the expression for $\mathbf{V_{22}}$ in~\eqref{eq:Vexp} and then enters the right hand side of the reconstruction formula of $\sigma$ in~\eqref{eq:grad-sigma}. We see from the figure that there appear small discontinuities at locations similar to the artifacts along the boundary $\partial \Omega \backslash \Gamma$. To obtain $\mathbf{V_{22}}$ one needs to compute the gradient of this function, which will cause the artifacts visible in the reconstruction. To support the theory that these are induced by the discretization of the domain, we compare the reconstruction of $\sigma$ for finer meshes in Figure~\ref{fig:mesh} and show the relative errors in Table~\ref{tab:RelErrMesh}. The mesh refinement does not make the artifacts disappear, but with increasing mesh size, the artifacts have less impact on the reconstruction, so that the relative errors decrease. 

\begin{table}[!ht]
    \centering
    \caption{Relative errors in the reconstructions for varying grid sizes $N_{\text{data}}$ and $N_{\text{recon}}$, using test case 1 and controlling $\Gamma_{\text{medium}}$.}
    \begin{tabular}{c|c|c|c}
        & \multicolumn{2}{c|}{number of nodes} & \\
        Grid size  & $N_{\text{data}}$ & $N_{\text{recon}}$ & Relative $L^2$ error $\sigma$ \\ \hline
        $N_{\text{small}}$ & 44.880 & 20.100 & 103.7\%\\
        $N_{\text{medium}}$ & 79.281 & 44.880 & 86.38 \%\\
        $N_{\text{large}}$ & 124.265 & 79.281 & 78.85\%\\\hline
    \end{tabular}
    \label{tab:RelErrMesh}
\end{table}

\subsection{Reconstruction of \texorpdfstring{$\sigma$}{sigma} from noisy data}
\begin{figure}[!ht]
    \centering
        \begin{minipage}[t]{0.33\textwidth}
        \centering
        \includegraphics[width=\linewidth]{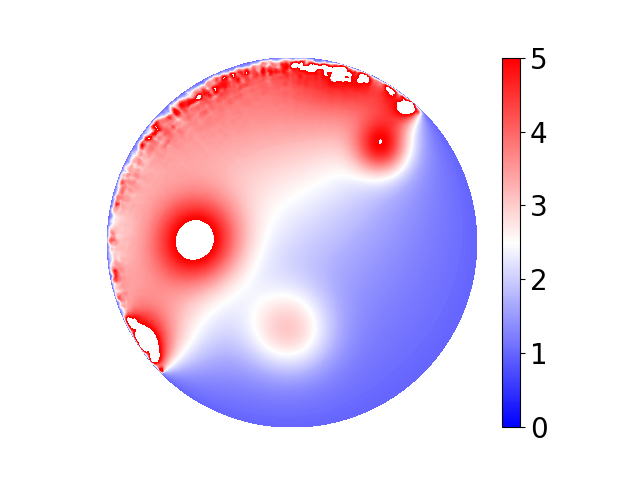}
        \caption*{$L=10^{-6}$}
    \end{minipage}%
    \begin{minipage}[t]{0.33\textwidth}
        \centering
        \includegraphics[width=\linewidth]{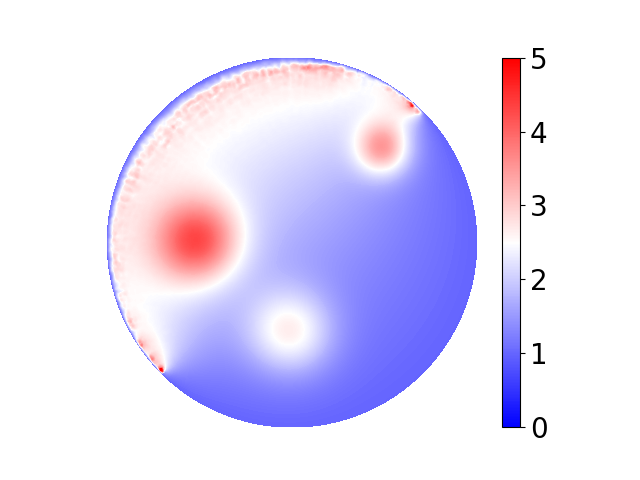}
        \caption*{$L=10^{-5}$}
    \end{minipage}%
    \begin{minipage}[t]{0.33\textwidth}
        \centering
        \includegraphics[width=\linewidth]{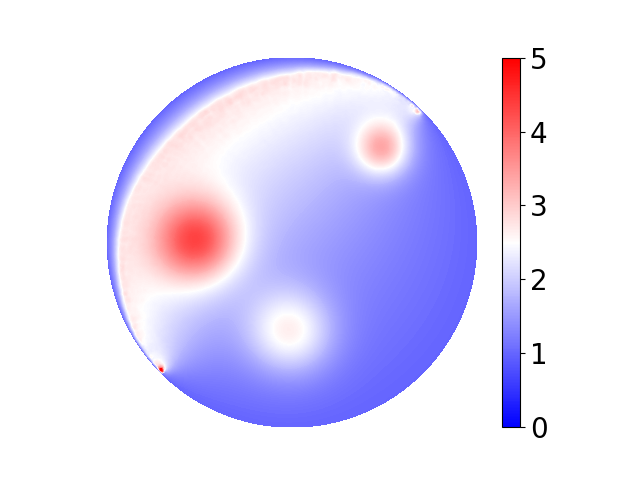}
        \caption*{$L=10^{-4}$}
    \end{minipage}
    \caption{Reconstructions of $\sigma$ as in test case 2 when corrupting the power density matrix with 5\% noise and using the boundary of control $\Gamma_{\text{medium}}$. We use varying lower bounds $L$ for the eigenvalues of $\widetilde{\mathbf H}$.} 
    \label{fig:RecNoiseL}
\end{figure}

\begin{figure}[!ht]
    \centering
        \begin{minipage}[t]{0.33\textwidth}
        \centering
        \includegraphics[width=\linewidth]{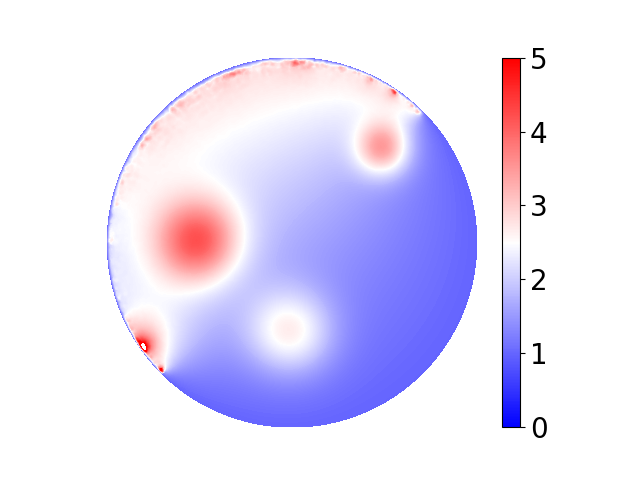}
        \caption*{1\% noise}
    \end{minipage}%
    \begin{minipage}[t]{0.33\textwidth}
        \centering
        \includegraphics[width=\linewidth]{FiguresNoise/RecSimga5noise1e5v2.png}
        \caption*{5\% noise}
    \end{minipage}%
    \begin{minipage}[t]{0.33\textwidth}
        \centering
        \includegraphics[width=\linewidth]{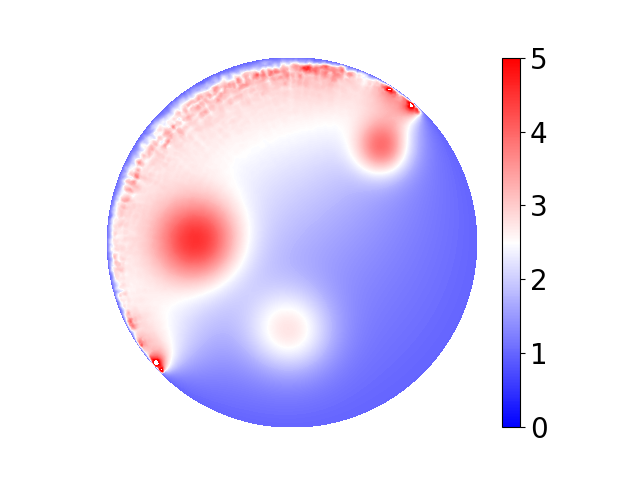}
        \caption*{10\% noise}
    \end{minipage}
    \caption{Reconstructions of $\sigma$ as in test case 2 when corrupting the power density matrix with varying noise levels and using the boundary of control $\Gamma_{\text{medium}}$. We use lower bounds $L$ for the eigenvalues of $\widetilde{\mathbf H}$: $L=10^{-6}$ (1\% noise) and $L=10^{-5}$ (5\% and 10\% noise).} 
    \label{fig:RecNoise}
\end{figure}

We perturb the entries of the power density matrix $\mathbf{H}$ at each node with random noise:
\begin{equation*}
    \widetilde{H}_{ij} = H_{ij} + \frac{\alpha}{100} \frac{e_{ij}}{\norm{e_{ij}}_{L^2}} H_{ij},
\end{equation*}
where $\alpha$ is the noise level and $e_{ij}$ are entries in the matrix $\mathbf{E}$ that are normally distributed $e_{ij} \sim \mathcal{N}(0,1)$. We use \texttt{numpy.random.randn} to generate the elements $e_{ij}$ and fix the seed \texttt{numpy.random.seed(50)}. The perturbation by random noise causes the Jacobian constraint to be violated close to the boundary $\partial \Omega \backslash \Gamma$. To enforce symmetry of $\mathbf{\widetilde{H}}$ we use $\frac{1}{2}(\mathbf{\widetilde{H}}+\mathbf{\widetilde{H}}^T)$ in the reconstruction procedure. Furthermore, to enforce positive definiteness of $\widetilde{\mathbf{H}}$ we use a small positive lower bound for the eigenvalues of $\widetilde{\mathbf{H}}$. This procedure works like a regularization of the reconstruction, where the lower bound, $L$, is the regularization parameter. This behavior is illustrated for reconstructions of $\sigma_{\text{Case 2}}$ in figure \ref{fig:RecNoiseL} and three different lower bounds $L=10^{-6}$, $L=10^{-5}$ and $L=10^{-4}$. Here $\mathbf H$ is perturbed with 5\% noise and $\Gamma_{\text{medium}}$ is used. When $L$ is chosen too small relative to the noise level the reconstruction is dominated by the noise on the values of $\widetilde{\mathbf{H}}$ for which $\det(\widetilde{\mathbf{H}})$ is small. On the other hand, for a large value of $L$ (like $L=10^{-4}$ in this case) a lot of information close to $\partial \Omega \backslash \Gamma$ is discarded, where the values of $\det(\widetilde{\mathbf{H}})$ should be small, as the Jacobian constraint is violated at the boundary. Therefore, we choose the lower bound as small as possible, while still obtaining a reasonable reconstruction which is the case for $L=10^{-5}$ in this case. Following this procedure, we reconstruct $\sigma_{\text{Case 2}}$ for three different noise levels $1\%, 5\%$ and $10\%$ in Figure~\ref{fig:RecNoise}. Here the lower bounds $L=10^{-6}$ (1\% noise) and $L=10^{-5}$ (5\% and 10\% noise) are used. Visually there appear more artifacts close to the boundary of control for increasing noise level. However, by choosing the suitable lower bound for each noise level the reconstructions are of similar quality. This makes it possible to add even higher levels of noise than 10\%, but for increasing noise level one needs to use increasing lower bounds.   

\section{Discussion and further work}
In this work we have shown the existence of finitely many boundary conditions that make solutions to the governing elliptic PDE satisfy the Jacobian constraint. This is certainly relevant for AET, as it indicates that AET in the limited view setting is indeed feasible. Moreover, we consider a numerical reconstruction framework in this setting and address both possibilities and limitations for reconstruction based on three particular boundary conditions.  

The results are based on the Runge approximation and hence they are not constructive. This immediately prompts the question of how to constructively find boundary conditions to satisfy the Jacobian constraint. An additional question concerns the reconstruction possibilities in the pragmatic choice of few boundary conditions. This is a question of deterioration of information with distance from the active part of the boundary as observed in~\cite{hubmer2018limited} and~\cite{cristo2019interior}. We note additionally, that given the no-flux condition outside $ \Gamma, $ the electrical field on the boundary must be tangential; hence the Jacobian constraint can never be satisfied on this part of the boundary. The reconstruction problem thus gets hard in the proximity to the boundary.

The result here itself guarantees a limited number of required boundary conditions, however, it is not clear how large that number may be. It would be interesting to find an upper bound at least for certain geometries.

A final important question for further studies concerns the growth and behavior of the boundary conditions found via the Runge approximation. The Runge approximation has recently been studied quantitatively~\cite{ruland2019quantitative}, and combining such results with our work could be an important step further.
    
\begin{acknowledgements}
BJ was supported by the Academy of Finland (grant no. 320022).
\end{acknowledgements}    
    
\clearpage
\bibliographystyle{mscplain}
\bibliography{bibliography}

\clearpage
\appendix

\section{Proof of Lemma~\ref{lem:appx-lem}}

\begin{proof}[Proof of Lemma~\ref{lem:appx-lem}]
    We sketch the principal steps of the proof. It is derived from~\cite[Section 7.3]{alberti2018lectures}.
    
    Step 1.) Observe first that for any $ g \in C^{0,\alpha}(B(x_0,s)) $ we have
    \begin{equation}
        \|g-g(x_0)\|_{C^{0,\frac\alpha2}(\overline{B(x_0,s)})} \leq C_1s^\frac{\alpha}{2}|g|_{C^{0,\alpha}(B(x_0,s))}, \label{eq:Holder-s-ineq}
    \end{equation}
    where $ C_1 $ is an independent constant.~\cite[Eq. (7.28)]{alberti2018lectures}
    
    Step 2.) Let $ u_s $ be the solution to the PDE problem
    \begin{subequations}
        \begin{alignat*}{2}
            Lu_s &= 0 &\quad& \text{in $ B(x_0,s) $}, \\ 
            u_s &= u_0 && \text{on $ \partial B(x_0,s) $}.
        \end{alignat*}
    \end{subequations}
    Consider now $ v_s = u_s - u_0 $ and observer that it is the unique solution of
    \begin{subequations}
        \begin{alignat*}{2}
            Lv_s &= -\nabla\cdot(\sigma -\sigma(x_0))\nabla u_0 &\quad& \text{in $ B(x_0,s) $}, \\ 
            v_s &= 0 && \text{on $ \partial B(x_0,s) $}.
        \end{alignat*}
    \end{subequations}
    Then by~\cite[Lem. 7.14]{alberti2018lectures} we have the estimate
    \begin{align*}
        \|v_s\|_{C^{1,\frac\alpha2}(\overline{B(x_0,s)})} \leq C_2\|(\sigma-\sigma(x_0))\nabla u_0\|_{C^{0,\frac\alpha2}(\overline{B(x_0,s)})}.
    \end{align*}
    By~\cite[Eq. (4.7)]{gilbarg2015elliptic} and~\eqref{eq:Holder-s-ineq}
    \begin{align*}
        \|(\sigma-\sigma(x_0))\nabla u_0\|_{C^{0,\frac\alpha2}(\overline{B(x_0,s)})} 
            \leq \|\sigma-\sigma(x_0)\|_{C^{0,\frac\alpha2}(\overline{B(x_0,s)})}\|\nabla u_0\|_{C^{0,\frac\alpha2}(\overline{B(x_0,s)})} 
            \leq C_3 s^{\frac\alpha2}
    \end{align*}
    where $ C_3 = C_1\|\nabla u_0\|_{C^{0,\frac\alpha2}(\overline{B(x_0,s)})}|\sigma|_{C^{0,\alpha}(B(x_0,s))} $. Therefore
    \begin{align*}
        \|u_s-u_0\|_{C^{1,\frac\alpha2}(\overline{B(x_0,s)})} \leq C_4 s^\frac{\alpha}{2}
    \end{align*}
    with $ C_4 = C_2C_3 $, and there is $ \tilde s \in (0,s) $ such that
    \begin{equation}
        \|u_{\tilde s}-u_0\|_{C^{1,\frac\alpha2}(\overline{B(x_0,\tilde{s})})} \leq \frac{\delta}{2}. \label{eq:appx-prf-delta-01}
    \end{equation}
    
    Step 3.) Observe that $ u_{\tilde s} \in H^1(B(x_0,\tilde{s})) \cap C^{1,\frac\alpha2}(B(x_0,\tilde{s})) $ and $ Lu_{\tilde{s}} = 0 $ in $ B(x_0,\tilde{s}) $. By Theorem~\ref{thm:L-runge-approx-weak} there is a sequence $ u_n \in H^1(\Omega) $ satisfying $ Lu_n = 0 $ in $ \Omega $ and $ L_\nu u_n = 0 $ on $ \partial\Omega\backslash \Gamma $ and $ u_n\vert_{B(x_0,\tilde{s})} \to u_{\tilde{s}} $ in $ L^2(B(x_0,\tilde{s})) $.
    
    Put $ w_n = u_n - u_{\tilde{s}} $ then $ w_n \in H^1(B(x_0,\tilde{s})) $ and satisfy $ Lw_n = 0 $ in $ B(x_0,\tilde{s}) $. By~\cite[Thm. 5.19]{giaquinta2013introduction} we have
    \begin{equation}
        \|\nabla w_n\|_{C^{0,\frac\alpha2}(\overline{B(x_0,\frac{\tilde{s}}{3})})} \leq C_5\|\nabla w_n\|_{L^2(B(x_0,\frac{\tilde{s}}{2})} \label{eq:appx-prf-ineq-01}
    \end{equation}
    with $ C_5 \equiv C_5(x_0,\tilde{s},\sigma) $, thus $ w_n \in C^{1,\frac\alpha2}(\overline{B(x_0,\frac{\tilde{s}}{3})}) $. By~\cite[Thm. 8.25]{salsa2016partial} we have 
    \begin{equation}
        \|w_n\|_{H^2(B(x_0,\frac{\tilde{s}}{k+1}))} \leq C_6\|w_n\|_{L^2(B(x_0,\frac{\tilde{s}}{k}))}, \label{eq:appx-prf-ineq-02}
    \end{equation}
    where $ C_6 \equiv C_6(x_0,\tilde{s},k,\sigma) $, and by Sobolev embedding~\cite{adams2003sobolev} we have
    \begin{equation}
        \|w_n\|_{C^0(B(x_0,\frac{\tilde{s}}{3}))} \leq C_7\|w_n\|_{H^2(B(x_0,\frac{\tilde{s}}{3}))} \label{eq:appx-prf-ineq-03}
    \end{equation}
    where $ C_7 \equiv C_7(x_0,\tilde{s}) $.
    
    By combining the estimates~\eqref{eq:appx-prf-ineq-01}--\eqref{eq:appx-prf-ineq-03} we obtain
    \[
        \|w_n\|_{C^{1,\frac\alpha2}(B(x_0,\frac{\tilde{s}}{3}))} \leq C_8\|w_n\|_{L^2(B(x_0,\tilde{s}))} \to 0 \quad \text{as $ n \to \infty $},
    \]
    where $ C_8 = C_8(x_0,\tilde{s},\sigma) $. Thus we have the strengthened convergence $ u_n \to u_{\tilde s} $ in $ C^{1,\frac\alpha2}(B(x_0,\frac{\tilde s}{3})) $. Hence we can pick $ \tilde n $ such that
    \begin{equation}
        \|u_{\tilde n} - u_{\tilde s}\|_{C^{1,\frac\alpha2}(B(x_0,\frac{\tilde{s}}{3}))} \leq \frac{\delta}{2}. \label{eq:appx-prf-delta-02}
    \end{equation}
    Step 4.) Choose $ r \in (0,\frac{\tilde{s}}3) $. Combining the estimates~\eqref{eq:appx-prf-delta-01} and~\eqref{eq:appx-prf-delta-02} we then have
    \[
        \|u_{\tilde{n}} - u_0\|_{C^1(\overline{B(x_0,r)})} \leq \|u_{\tilde{n}} - u_{\tilde s}\|_{C^{1,\frac\alpha2}(\overline{B(x_0,r)})} + \|u_{\tilde{s}} - u_0\|_{C^{1,\frac\alpha2}(\overline{B(x_0,r)})} \leq \frac{\delta}{2} + \frac{\delta}{2} = \delta.
    \]
    
\end{proof} 

\section{Weak Maximum Principle for a mixed Dirichlet-Neumann problem}
\begin{theorem}[Weak Maximum Principle]\label{thm:Max}
    Let $\Omega \subset \mathbb{R}^n$ be a bounded connected Lipschitz domain and $\sigma \in L^{\infty}_{+}(\Omega)$. Let $u \in H^1(\Omega) \cap C(\overline{\Omega})$ be a solution to \eqref{eq:u} so that it satisfies the variational equation
    \begin{equation}\label{eq:vareq}
        \int_{\Omega} \sigma \nabla u \cdot \nabla v \, \mathrm{d}x = 0,
    \end{equation}
    for all test functions $v \in H_{0,\Gamma}^1(\Omega)$. Here $H_{0,\Gamma}^1(\Omega)$ consists of the set of functions in $H^1(\Omega)$ with zero trace on $\Gamma$. Then the following are true:
    \begin{enumerate}
        \item If $u \geq 0$ on $\Gamma$ then $u \geq 0$ in $\Omega$.
        \item $\inf_{y \in \Gamma} u(y) \leq u(x) \leq \sup_{y \in \Gamma} u(y)$ for all $x \in \Omega$.
    \end{enumerate}
\end{theorem}
\begin{proof}
    This proof is based on \cite[proof of Thm. 8.20]{salsa2016partial}, but for the convenience of the reader it is simplified to only cover the setting where the elliptic operator is composed of a diffusion term. Define $u^-=\max\{-u,0\}$ and $u^+=\max\{u,0\}$. Since $u \in H^1(\Omega)$ then $\vert u \vert$, $u^-$ and $u^+$ are also in $H^1(\Omega)$ \cite[Prop. 7.68]{salsa2016partial}. 
    \begin{enumerate}
        \item Since $u \geq 0$ on $\Gamma$, the negative part $u^-$ has zero trace on $\Gamma$. Hence, $u^- \in H_{0,\Gamma}^1$. Thus for the choice of the test function $v=u^-$, the variational equation \eqref{eq:vareq} is 
        \begin{align*}
            \int_{\Omega} \sigma \nabla u \cdot \nabla u^- \, \mathrm{d}x = 0.
            \intertext{This simplifies to the following as $\nabla u^- =\nabla u \,\chi_{u<0}$ a.e. in $\Omega$, so there is only a contribution for $u<0$ and in this case $u=-u^-$:}
            -\int_{\Omega} \sigma \nabla u^- \cdot \nabla u^- \, \mathrm{d}x = 0.
        \end{align*}
        As $\sigma \in L^{\infty}_+(\Omega)$ it follows directly from this equation that $\nabla u^-=0$ in $L^2(\Omega)$, implying $u^-=C$, where $C$ is a constant. As $u^-$ has zero trace on $\Gamma$ it follows that $C=0$. This implies $u^-=0$ a.e. and therefore $u=u^+$ a.e. on $\Omega$. 
        
        \item Let $m_1=\inf_{y \in \Gamma} u(y)$. Then $u-m_1 \geq 0$ and $u-m_1 \in H^1(\Omega)$ satisfies the variational equation \eqref{eq:vareq}, as $m_1$ is just a constant. By step (1) it then follows that $u-m_1 = (u-m_1)^+$ a.e. on $\Omega$, hence $u \geq \inf_{y \in \Gamma} u(y)$ a.e. on $\Omega$. Defining $m_2=\sup_{y \in \Gamma} u(y)$ and using the same approach for the function $m_2-u$ yields $u \leq \sup_{y \in \Gamma} u(y)$ a.e. on $\Omega$.
    \end{enumerate}
\end{proof}

\end{document}